# ON THE POWER OF TWO CHOICES: BALLS AND BINS IN CONTINUOUS TIME

By Malwina J. Luczak and Colin McDiarmid

*London School of Economics and University of Oxford*

Suppose that there are $n$ bins, and balls arrive in a Poisson process at rate $\lambda n$, where $\lambda > 0$ is a constant. Upon arrival, each ball chooses a fixed number $d$ of random bins, and is placed into one with least load. Balls have independent exponential lifetimes with unit mean. We show that the system converges rapidly to its equilibrium distribution; and when $d \geq 2$, there is an integer-valued function $m_d(n) = \ln\ln n/\ln d + O(1)$ such that, in the equilibrium distribution, the maximum load of a bin is concentrated on the two values $m_d(n)$ and $m_d(n) - 1$, with probability tending to 1, as $n \to \infty$. We show also that the maximum load usually does not vary by more than a constant amount from $\ln\ln n/\ln d$, even over quite long periods of time.

**1. Introduction.** Balls-and-bins processes have been useful for modeling and analyzing a wide range of problems, in discrete mathematics, computer science and communication theory, and, in particular, for problems which involve load sharing, see, for example, [4, 5, 12, 15, 16, 17, 22]. Here is one central result, from [3]. Let $d$ be a fixed integer at least 2. Suppose that there are $n$ bins, and $n$ balls arrive one after another: each ball picks $d$ bins uniformly at random and is placed in a least loaded of these bins. Then with probability tending to 1 as $n \to \infty$, the maximum load of a bin is $\ln\ln n/\ln d + O(1)$.

In some recent work, balls have been allowed to "die," see [3, 7, 21], which is, of course, desirable when modeling telephone calls. For example, suppose that we start with $n$ balls in $n$ bins: at each time step, one ball is deleted uniformly at random, and one new ball appears and is placed in one of $d$ bins as before. It is shown in [3] that, as $n \to \infty$, at any given time $t \geq cn^2 \ln \ln n$, with probability tending to 1, the maximum load of a bin is at most $\ln \ln n / \ln d + O(1)$.

The results mentioned above all concern discrete time models, where at each time step a ball may arrive or a ball may die and be replaced by a new one. Here we analyze a simple and natural continuous time "immigration–death" balls-and-bins model. We concentrate on the maximum bin load, which may be the quantity of greatest interest, for example, in load-sharing models.

The scenario we consider is as follows. Let $d$ be a fixed positive integer, say $d = 2$. Let $n$ be a positive integer and suppose that there are $n$ bins. Balls arrive in a Poisson process at rate $\lambda n$, where $\lambda > 0$ is a constant. Upon arrival, each ball chooses $d$ random bins (with replacement), and is placed into a least loaded bin among those chosen. (If there is more than one chosen bin with least load, the ball is placed in the first such bin chosen.) Balls have independent exponential lifetimes with unit mean. This process goes on forever.

This model was first studied by Turner in [21], who considers weak convergence, for a suitable choice of state space. (Also, [19, 20] contain a discussion of the completeness of the state space under the product topology.) Turner shows that (with appropriate assumptions on the initial distribution), for each fixed non-negative integer $k$ the fraction of bins with load at least $k$, converges weakly as $n \to \infty$ to a deterministic function $v(t,k)$ defined on $\mathbb{R}^+ \times \mathbb{Z}^+$, where the vector $(v(t,k) : k \in \mathbb{Z}^+)$ is the unique solution to the system of differential equations for $k = 1, 2, \ldots$,

$$(1) \quad \frac{dv(t,k)}{dt} = \lambda(v(t,k-1)^d - v(t,k)^d) - k(v(t,k) - v(t,k+1)), \qquad t \geq 0,$$

subject to $v(t,0) = 1$ for all $t \geq 0$, and appropriate initial values $(v(0,k) : k \in \mathbb{Z}^+)$ such that $1 \geq v(0,k) \geq v(0,k+1) \geq 0$ for all $k \in \mathbb{N}$. The weak-convergence result applies only to fixed-index co-ordinates (i.e., fixed values of $k$) over fixed-length time intervals, and yields no information on the speed of convergence. Our approach is different, and we are not concerned with weak convergence, although weak convergence could be deduced from our results. The key step is to establish concentration results, which apply to the fraction of bins with load at least $k$ at time $t$ (where $k, t$ need not be fixed); these concentration results may then be used to analyze a balance equation involving these quantities. We are thus able to handle random variables like the maximum load, over long periods of time.

For each time $t \geq 0$ and each $j = 1, \ldots, n$, let $X_t(j)$ be the random number of balls in bin $j$ at time $t$, and let $X_t$ be the *load vector* $(X_t(1), \ldots, X_t(n))$. Thus, the total number of balls $|X_t|$ at time $t$ is given by $|X_t| = \sum_{j=1}^n X_t(j)$. We shall always assume that the initial load vector $X_0$ satisfies $\mathbf{E}[|X_0|] <$

$\infty$. Note that $|X_t|$ follows a simple immigration–death process, and so its stationary distribution is the Poisson distribution $\text{Po}(\lambda n)$ with mean $\lambda n$.

It is easy to check that, for given $d$ and $n$, the load vector process $(X_t)$ is Markov, with state space $(\mathbb{Z}^+)^n$. Standard results show that there is a unique stationary distribution $\mathbf{\Pi}$; and, whatever the distribution of the starting state $X_0$, the distribution of the load vector $X_t$ at time $t$ converges to $\mathbf{\Pi}$ as $t \to \infty$. Indeed, this convergence is very fast, as our first theorem will show.

For $x \in \mathbb{Z}^n$, let $\|x\|_1 = \sum_i |x(i)|$ be the $L_1$ norm of $x$. (Thus, we have $|X_t| = \|X_t\|_1$.) We use $\mathcal{L}(X)$ to denote the probability law or distribution of a random variable $X$. The *total variation distance* between two probability distributions $\mu_1$ and $\mu_2$ may be defined by $d_{\text{TV}}(\mu_1, \mu_2) = \inf \mathbf{Pr}(X \neq Y)$, where the infimum is over all couplings of $X$ and $Y$, where $\mathcal{L}(X) = \mu_1$ and $\mathcal{L}(Y) = \mu_2$. Equivalently,

$$d_{\text{TV}}(\mu_1, \mu_2) = \max_A |\mathbf{Pr}(X \in A) - \mathbf{Pr}(Y \in A)|,$$

where the maximum is over all suitable sets $A$. We also use the Wasserstein distance, defined by $d_{\text{W}}(\mu_1, \mu_2) = \inf \mathbf{E}[\|X - Y\|_1]$, where the inf is over couplings of $X$ and $Y$ as above. For distributions $\mu_1$ and $\mu_2$ on $\mathbb{Z}^n$, we have $d_{\text{TV}}(\mu_1, \mu_2) \leq d_{\text{W}}(\mu_1, \mu_2)$.

THEOREM 1.1. *Let $d$ and $n$ be positive integers, and let $\mathbf{\Pi}$ be the corresponding stationary distribution for the load vector. Suppose that initially the balls are arbitrarily distributed over the bins, with $\mathbf{E}[|X_0|] < \infty$. Then for each time $t \geq 0$,*

$$d_{\text{TV}}(\mathcal{L}(X_t), \mathbf{\Pi}) \leq d_{\text{W}}(\mathcal{L}(X_t), \mathbf{\Pi}) \leq (\lambda n + \mathbf{E}[|X_0|])e^{-t}.$$

For each $\varepsilon > 0$ and initial state $x$, the *mixing time* $\tau(\varepsilon, x)$ is defined by considering $(X_t)$, where $X_0 = x$ a.s. and setting

$$\tau(\varepsilon, x) = \inf\{t \geq 0 : d_{\text{TV}}(\mathcal{L}(X_t), \mathbf{\Pi}) \leq \varepsilon\}.$$

[Recall that $d_{\text{TV}}(\mathcal{L}(X_t), \mathbf{\Pi})$ is a nonincreasing function of $t$.] Thus, for example, if $\mathbf{0}$ denotes the state with no balls, then the above theorem shows that

$$\tau(\varepsilon, \mathbf{0}) \leq \ln(\lambda n / \varepsilon).$$

This upper bound on the mixing time is, in fact, of the right order, in that $\tau(\frac{1}{2}, \mathbf{0}) = \Theta(\ln n)$, as we shall see after the proof of Theorem 1.1 by considering the behavior of the total number of balls present. For mixing results on related models, see [4, 7]: mixing appears to be slower when balls live forever.

As we commented earlier, our primary interest is in the maximum load of a bin. Let $M_t = \max_j X_t(j)$ be the maximum load of a bin at time $t$.

Thus, $M_t = \|X_t\|_\infty$, where $\|x\|_\infty$ is the infinity norm $\max_j |x_j|$ of $x$. The above theorem shows that we can essentially restrict our attention to the stationary case, at least if we are interested in times well beyond $\ln n$, so let us now consider that case. We may write $M$ instead of $M_t$ when the system is in equilibrium. The behavior of the maximum load $M_t$ or $M$ is very different in the two cases $d = 1$ and $d \geq 2$. This is the "power of two choices" phenomenon—see, for example, [17]. For clarity, let us write $X_t^{(n)}$ and $M_t^{(n)}$ or $M^{(n)}$ here to indicate that there are $n$ bins.

The most interesting case is when $d \geq 2$ (indeed, when $d = 2$), but in order to set things in context, let us first consider the (much easier) case when $d = 1$. Suppose then that $d = 1$. We shall see that $M^{(n)}$ is concentrated on two values $m = m(n)$ and $m - 1$, which are close to $\ln n / \ln \ln n$; and that over a polynomial length interval of time, we meet only small (constant size) deviations below $m$ but we meet large deviations above $m$, so that the maximum value of $M_t^{(n)}$ over an interval of length $n^K$ is usually about $(K+1)m$. We use the phrase *asymptotically almost surely* (a.a.s.) to mean "with probability $\to 1$ as $n \to \infty$."

THEOREM 1.2. *Let $d = 1$, and suppose that $X_0^{(n)}$ is in the stationary distribution (and thus so is $M_t^{(n)}$ for each time $t$).*

(a) *There exists an integer-valued function $m = m(n) \sim \frac{\ln n}{\ln \ln n}$ such that a.a.s. $M^{(n)}$ is $m(n)$ or $m(n) - 1$.*

(b) *For any constant $K > 0$,*

$$\min_{0 \leq t \leq n^K} M_t^{(n)} \geq m(n) - 3 \qquad a.a.s.$$

(c) *For any constant $K > 0$,*

$$\left(\max_{0 \leq t \leq n^K} M_t^{(n)}\right) \frac{\ln \ln n}{\ln n} \to K + 1 \qquad \text{in probability as } n \to \infty.$$

The notation $m = m(n) \sim \frac{\ln n}{\ln \ln n}$ above means that $m(n) = (1 + o(1))\frac{\ln n}{\ln \ln n}$ as $n \to \infty$. It is straightforward to determine $m(n)$ more precisely from the proof of the theorem: for example, we have

$$m(n) = \frac{\ln n}{\ln \ln n} + (1 + o(1))\frac{(\ln n)(\ln \ln \ln n)}{(\ln \ln n)^2}.$$

Now we consider the case $d \geq 2$, when the maximum load $M_t^{(n)}$ is far smaller. Once again, it is concentrated on two values $m_d = m_d(n)$ and $m_d - 1$, but now these numbers are close to $\ln \ln n / \ln d$. This corresponds to the behavior of the maximum load in discrete time models; see, for example, [3, 4, 12, 16], but is more precise.

THEOREM 1.3. *Let $d \geq 2$ be fixed, and suppose that $X_0^{(n)}$ is in the stationary distribution. Then there exists an integer-valued function $m_d = m_d(n) = \ln \ln n / \ln d + O(1)$ such that $M^{(n)}$ is $m_d$ or $m_d - 1$ a.a.s. Further, for any constant $K > 0$, there exists $c = c(K)$ such that*

$$\max_{0 \leq t \leq n^K} |M_t^{(n)} - \ln \ln n / \ln d| \leq c \qquad a.a.s. \tag{2}$$

The lower bound on $M_t^{(n)}$, in fact, holds over longer intervals than stated in (2) above. For example, there is a constant $c$ such that

$$\min\{M_t^{(n)} : 0 \leq t \leq e^{n^{1/4}}\} \geq \ln \ln n / \ln d - c \qquad \text{a.a.s.} \tag{3}$$

However, the upper bound in (2) does not extend to much longer intervals. For example, if $K > 0$ and $\tau = n^{Kd \ln \ln n}$, then

$$\max_{0 \leq t \leq \tau} M_t^{(n)} \geq K \ln \ln n \qquad \text{a.a.s.} \tag{4}$$

The plan of the rest of the paper is as follows. After giving some preliminary results in the next section, we consider mixing times and prove Theorem 1.1. Then we consider the easy case $d = 1$ when there is one random choice, and prove Theorem 1.2. In order to prove Theorem 1.3, where $d \geq 2$, we need some preliminary results, which are presented in the next three sections. First, in Section 5 we give a concentration result for Lipschitz functions of the load vector in equilibrium. In Section 6 we use balance equations to establish the key equation (26) concerning the expected proportion $u(i)$ of bins with load at least $i$ in equilibrium. This result, together with the concentration result, yields a recurrence for $u(i)$. After that, in Section 7 we consider random processes like a random walk with "drift." Then we are ready to prove Theorem 1.3 in Section 8: we first prove upper bounds, then lower bounds, and finally we prove the results (3) and (4). Last, we briefly consider chaoticity and make some concluding remarks.

**2. Preliminary results.** In this section we give some elementary results which we shall need several times below. A standard inequality for a binomial or Poisson random variable $X$ with mean $\mu$ is that

$$\mathbf{Pr}(|X - \mu| \geq \epsilon \mu) \leq 2 \exp(-\tfrac{1}{3}\epsilon^2 \mu) \tag{5}$$

for $0 \leq \epsilon \leq 1$ (see, e.g., Theorem 2.3(c) and inequality (2.8) in [14]). Also, for each positive integer $k$,

$$\mathbf{Pr}(X \geq k) \leq \mu^k / k! \leq (e\mu/k)^k. \tag{6}$$

If $X$ has the Poisson distribution with mean $\mu$, let us write $X \sim \text{Po}(\mu)$: for such a random variable, we have

$$\mathbf{E}[X \mathbf{1}_{(X \geq k)}] = \mu \mathbf{Pr}(X \geq k - 1). \tag{7}$$

Next we give an elementary lemma which we shall use later in order to extend certain results, for example, concerning the maximum load $M_t$ from a single point in time to an interval of time. It yields bounds on the maximum and minimum values of a suitable function $f(x)$ over a time interval $[0, \tau]$.

Consider the $n$-bin case, with set $\Omega = (\mathbb{Z}^+)^n$ of load vectors. Let us say that a real-valued function $f$ on $\Omega$ has *bounded increase* if whenever $s$ and $t$ are times with $s < t$, then $f(x_t)$ is at most $f(x_s)$ plus the total number of arrivals in the interval $(s,t]$; and $f$ has *strongly bounded increase* if $f(x_t)$ is at most $f(x_s)$ plus the maximum number of arrivals in the interval $(s,t]$ which are placed in any one bin. Thus, for example, $f(x) = |x|$ has bounded increase, and $f(x) = \max_j x(j)$ has strongly bounded increase.

LEMMA 2.1. *Let $(X_t)$ be in equilibrium. Let $s, \tau > 0$ and let $a, b$ be non-negative integers. Suppose that* (a) *$f$ has bounded increase and $\delta = \mathbf{Pr}(\mathrm{Po}(\lambda n s) \geq b+1)$, or* (b) *$f$ has strongly bounded increase and $\delta = n \mathbf{Pr}(\mathrm{Po}(\lambda ds) \geq b+1)$. In both cases we have*

$$(8) \quad \mathbf{Pr}\left[f(X_t) \leq a \text{ for some } t \in [0, \tau]\right] \leq \left(\frac{\tau}{s} + 1\right)\left(\mathbf{Pr}\left(f(X_0) \leq a + b\right) + \delta\right)$$

*and*

$$(9) \quad \mathbf{Pr}\left[f(X_t) \geq a + b \text{ for some } t \in [0, \tau]\right] \leq \left(\frac{\tau}{s} + 1\right)\left(\mathbf{Pr}\left(f(X_0) \geq a\right) + \delta\right).$$

PROOF. Consider first the case (a) when $f$ has bounded increase. Note that the $j = \lfloor \frac{\tau}{s} \rfloor + 1$ disjoint intervals $[(r-1)s, rs)$ for $r = 1, \ldots, j$ cover $[0, \tau]$. Let $B_r$ denote the event of having in total at least $b+1$ arrivals in the interval $[(r-1)s, rs)$, so that $\mathbf{Pr}(B_r) = \mathbf{Pr}[\mathrm{Po}(\lambda n s) \geq b+1] = \delta$. Then

$$\{f(X_t) \leq a \text{ for some } t \in [0, \tau]\} \subseteq \left(\bigcup_{r=1}^{j} \{f(X_{rs}) \leq a + b\}\right) \cup \left(\bigcup_{r=1}^{j} B_r\right)$$

and (8) follows. Similarly,

$$\{f(X_t) \geq a + b \text{ for some } t \in [0, \tau)\} \subseteq \left(\bigcup_{r=0}^{j-1} \{f(X_{rs}) \geq a\}\right) \cup \left(\bigcup_{r=1}^{j} B_r\right)$$

and (9) follows. To handle the case (b) when $f$ has strongly bounded increase, note that if $C_r$ denotes the event of having at least $b+1$ arrivals in the interval $[(r-1)s, rs)$ which are placed into a single bin, then $\mathbf{Pr}(C_r) \leq n \mathbf{Pr}[\mathrm{Po}(\lambda ds) \geq b+1]$; and then proceed as above. $\square$

As we noted earlier, in equilibrium the distribution of the total number of balls in the system is $\mathrm{Po}(\lambda n)$. We close this section by using the last lemma to establish a result that will enable us to "control" the total number of balls in the system over long periods of time.

LEMMA 2.2. *For any $0 < \epsilon < 1$, there exists $\beta > 0$ such that the following holds. Consider an $n$-bin system, and let $(X_t)$ be in equilibrium. Then a.a.s. for all $0 \leq t \leq e^{\beta n}$, the number of balls $|X_t|$ satisfies*

$$(1-\epsilon)\lambda n \leq |X_t| \leq (1+\epsilon)\lambda n.$$

PROOF. By inequality (5), since $|X_t| \sim \text{Po}(\lambda n)$, we have

$$\mathbf{Pr}\left(||X_t| - \lambda n| > \epsilon \lambda n/2\right) \leq 2e^{-\epsilon^2 \lambda n/12}$$

and

$$\mathbf{Pr}\left[\text{Po}(\epsilon \lambda n/4) \geq \epsilon \lambda n/2\right] \leq 2e^{-\epsilon \lambda n/12}.$$

Let $\beta$ satisfy $0 < \beta < \frac{1}{12}\epsilon^2 \lambda$. We use case (a) of Lemma 2.1. Let $s = \epsilon/4$ and $b = \epsilon \lambda n/2$: we may now use (8) with $a = (1-\epsilon)\lambda n$ and (9) with $a = (1+\epsilon/2)\lambda n$.  □

**3. Rapid mixing: proof of Theorem 1.1.** We shall couple $(X_t)$ and a corresponding copy $(Y_t)$ of the process in equilibrium in such a way that with high probability $\|X_t - Y_t\|_1$ decreases quickly to 0. We assume that the choices process always generates a nonempty list of bins at an arrival time, and the new ball is placed in a least-loaded bin among those chosen, breaking ties if necessary by choosing the first least-loaded bin in the list. In the meantime we make no other assumptions about the arrivals process or the choices process. We assume as before that balls die independently at rate 1, independently of the other two processes.

The coupling is as follows. Not surprisingly, we give the two processes the same arrivals and choices of $d$ bins. The *height* of a ball in the system at a given time is the number of balls in its bin that arrived before it, plus one. Assume that we have a family of independent rate 1 Poisson processes $F_{j,k}$ for $j = 1, \ldots, n$ and $k = 1, 2, \ldots$. When $F_{j,k}$ "tolls," any ball in bin $j$ at height $k$ in either process dies (so that 0 or 1 or 2 balls die). Observe that at any time $t$, we are interested in only a finite (with probability 1) number of these death processes [namely, $\sum_j X_t(j) \vee Y_t(j)$]. We have now described the coupling of $(X_t)$ and $(Y_t)$. The "memoryless" property of the exponential lifetime distribution ensures that it is a proper coupling; and when the arrival process is Poisson, and the choices are independent and uniform, the joint process $(X_t, Y_t)$ is Markov. For $x, y \in \mathbb{Z}^n$, the notation $x \leq y$ means that $x(j) \leq y(j)$ for each $j = 1, \ldots, n$.

LEMMA 3.1. *With the coupling of $(X_t)$ and $(Y_t)$ described above, the distance $\|X_t - Y_t\|_1$ is nonincreasing, and given that $\|X_0 - Y_0\|_1 = r$, it is stochastically at most the number of survivors at time $t$ of $r$ independent balls. Further, if $0 \leq s \leq t$ and $X_s \leq Y_s$, then $X_t \leq Y_t$.*

PROOF. Consider a jump time $t_0$. Let $X_{t_0-} = x$ and $Y_{t_0-} = y$, and let $X_{t_0} = x'$ and $Y_{t_0} = y'$. (We assume right-continuity.) Suppose that $t_0$ is a death ("toll") time. If none or two balls die, then

$$\|x' - y'\|_1 = \|x - y\|_1, \tag{10}$$

and if just one ball dies, then

$$\|x' - y'\|_1 = \|x - y\|_1 - 1. \tag{11}$$

Thus, at any death time,

$$\|x' - y'\|_1 \leq \|x - y\|_1. \tag{12}$$

Suppose now that $t_0$ is an arrival time, and ball $b$ arrives. We want to show that (12) holds. If ball $b$ is placed in the same bin in the two processes, then (10) holds and, hence, so does (12). Suppose that ball $b$ is placed in bin $i$ in the $X$-process and in bin $j$ in the $Y$-process, where $i \neq j$. Then ball $b$ gets "paired" in at least one of the processes, and so (12) holds. (By "paired" here, we mean that in the other process there is a ball in the same bin at the same height. Observe that these balls will stay paired until they die together.) For, note first that $x(i) \leq x(j)$ and $y(j) \leq y(i)$, and not both are equal by the tie-breaking rule. Now suppose that ball $b$ does not get paired in either process. Then we must have $x(i) \geq y(i)$ and $y(j) \geq x(j)$, and so

$$x(i) \geq y(i) \geq y(j) \geq x(j) \geq x(i).$$

But then all the values are equal, a contradiction.

We have now seen that (12) holds at each jump time, and (11) holds if a single unpaired ball dies. Thus, $\|X_t - Y_t\|_1$ is nonincreasing. Further, we claim that, for any time $0 \leq s < t$ and any positive integer $r$, given that $\|X_s - Y_s\|_1 = r$ and any other history up to time $s$, the probability that $\|X_t - Y_t\|_1 = r$ is at most $e^{-r(t-s)}$. The second part of the lemma will follow immediately from the claim.

To see why the claim is true, let $S_r$ denote the set of states $(x, y)$ such that $\|x - y\|_1 = r$. We have seen that $\|X_t - Y_t\|_1$ is nonincreasing. For each state $(x, y) \in S_r$, there are $r$ of the death processes $F_{jk}$ such that if any of them tolls, then the process moves into $S_{r-1}$. Thus, if $(X_0, Y_0) \in S_r$ and $T = \inf\{t \geq 0 : (X_t, Y_t) \notin S_r\}$ is the exit time from $S_r$, then $\mathbf{Pr}(T > t | (X_0, Y_0) = (x, y)) \leq e^{-rt}$ for each $(x, y) \in S_r$ and each $t > 0$; and the claim follows.

The final comment on monotonicity is straightforward. For consider a jump time $t_0$ as above, and suppose that $x \leq y$. If $t_0$ is a death time, then clearly $x' \leq y'$, so suppose that $t_0$ is an arrival time. But if the new ball is placed in bin $i$ in the $X$-process and if $x(i) = y(i)$, then the ball is placed in bin $i$ also in the $Y$-process, so $x' \leq y'$. □

We may now rapidly prove Theorem 1.1. By the lemma, $\mathbf{E}(\|X_t - Y_t\|_1 | (X_0, Y_0) = (x,y))$ is at most the expected number among $r = \|x - y\|_1$ balls that survive at least to time $t$, which is equal to $re^{-t}$. Since $\|x - y\|_1 \leq |x| + |y|$, we have

$$\mathbf{E}(\|X_t - Y_t\|_1 | X_0, Y_0) \leq (|X_0| + |Y_0|)e^{-t},$$

and so

$$d_W(\mathcal{L}(X_t), \mathcal{L}(Y_t)) \leq \mathbf{E}(\|X_t - Y_t\|_1) \leq (\mathbf{E}[|X_0|] + \lambda n)e^{-t}.$$

This completes the proof of Theorem 1.1.

We now show that the upper bounds on the mixing times arising from Theorem 1.1 are of the right order. We may see this by simply considering the total number $|X_t|$ of balls in the system. In equilibrium, $|X_t|$ has the Poisson distribution $\text{Po}(\lambda n)$, and so

$$d_{\text{TV}}(\mathcal{L}(X_t), \mathbf{\Pi}) \geq d_{\text{TV}}(\mathcal{L}(|X_t|), \text{Po}(\lambda n)).$$

We shall see that if $X_0 = \mathbf{0}$ a.s. and $t \leq \frac{1}{2}\ln n - 2\ln\ln n$, then

(13) $$d_{\text{TV}}(\mathcal{L}(|X_t|), \text{Po}(\lambda n)) = 1 - o(1);$$

and it follows that, for each $0 < \varepsilon < 1$, we have $\tau(\varepsilon, \mathbf{0}) = \Omega(\ln n)$.

Suppose then that $X_0 = \mathbf{0}$ a.s. and let $\mu(t) = \mathbf{E}[|X_t|]$. It is easy to check that $\mu(t) = \lambda n(1 - e^{-t})$. If $t$ is $\Theta(\ln n)$, then, by Lemma 5.5 below (with, say, $b = \ln^{3/2} n$),

$$\mathbf{Pr}\left(||X_t| - \mu(t)| \geq \tfrac{1}{2}\lambda n^{1/2}\ln^2 n\right) = e^{-\Omega(\ln^{3/2} n)}.$$

Also, if $Z \sim \text{Po}(\lambda n)$, then, by (5),

$$\mathbf{Pr}\left(|Z - \lambda n| \geq n^{1/2}\ln n\right) = e^{-\Omega(\ln^2 n)}.$$

Now if $t$ is $\frac{1}{2}\ln n - 2\ln\ln n$, then $|\mu(t) - \lambda n| = \lambda n e^{-t} = \lambda n^{1/2}\ln^2 n$, and, thus,

$$d_{\text{TV}}(\mathcal{L}(|X_t|), \text{Po}(\lambda n)) = 1 - e^{-\Omega(\ln^{3/2} n)} = 1 - o(1),$$

which gives (13) as required (since the left-hand side is a nonincreasing function of $t$).

**4. One choice: proof of Theorem 1.2.** Let $\lambda > 0$ be fixed, as always. Let $d = 1$. Let $p_i = p_i(\lambda) = e^{-\lambda}\sum_{k \geq i}\frac{\lambda^k}{k!}$, the probability that a $\text{Po}(\lambda)$ random variable takes value at least $i$. Let $X_0$ be in equilibrium. Stationary bin loads are independent Poisson random variables, each with mean $\lambda$. It follows that, for any nonnegative integer $i$,

(14) $$\mathbf{Pr}(M_t \geq i) \leq np_i$$

and

(15) $$\mathbf{Pr}(M_t \leq i) = (1 - p_{i+1})^n \leq e^{-np_{i+1}}.$$

We now prove the three parts of the theorem.

*Part* (a). Let $\omega(n) = \ln \ln n$. Let $m = m(n)$ be the least positive integer $i$ such that $np_{i+1} \leq 1/\omega(n)$. By (14),

$$\mathbf{Pr}\,(M_t \geq m+1) \leq np_{m+1} = o(1),$$

so $M_t \leq m$ a.a.s. Also, $np_m > 1/\omega(n)$, so $np_{m-1} = \Omega(\frac{\ln n}{\ln \ln n} \cdot \frac{1}{\omega(n)}) \to \infty$. Hence, by (15),

$$\mathbf{Pr}\,(M_t \leq m-2) \leq e^{-np_{m-1}} = o(1).$$

Thus, $M_t$ is $m$ or $m-1$ a.a.s. Also, it is easy to check that $m \sim \frac{\ln n}{\ln \ln n}$.

*Part* (b). We apply case (b) of Lemma 2.1, with $s \sim n^{-K-2}$, $a = m-4$ and $b = 1$, together with (6) and (15).

*Part* (c). Let $Z = \max_{0 \leq t \leq n^K} M_t$. Let $\varepsilon > 0$. We show first that

(16) $\quad \mathbf{Pr}\,(Z > (K+1+\varepsilon)\ln n/\ln \ln n) \to 0 \quad \text{as } n \to \infty.$

To do this, we apply case (b) of Lemma 2.1, with $s \sim \exp(-\ln n/\ln \ln n)$, $a \sim (K+1+\varepsilon/2)\ln n/\ln \ln n$ and $b \sim \ln n/(\ln \ln n)^2$, together with (6) and (14).

Now let $0 < \varepsilon < K$, and let $k = \lceil (K+1-\varepsilon)\ln n/\ln \ln n \rceil$. We will show that

(17) $\quad \mathbf{Pr}\,(Z < k) \to 0 \quad \text{as } n \to \infty,$

which will complete the proof of this part and thus of the theorem. Note that $np_k = n^{-(K-\varepsilon+o(1))} = o(1)$. For each time $t > 0$, let $\phi_t$ be the sigma field generated by all events until time $t$. Let $C$ be the event that $|X_t| \leq n^2/2$ for each $t \in [0, n^K]$. Then $C$ holds a.a.s. by Lemma 2.2. Let $n \geq 2\lambda$ and let $x$ be a load vector such that $|x| \leq n^2/2$. Given $X_0 = x$, by Theorem 1.1,

$$d_{\mathrm{TV}}(\mathcal{L}(X_t), \Pi) \leq (\lambda n + |x|)e^{-t}$$
$$\leq n^2 e^{-t}$$
$$\leq e^{-\ln^2 n}$$

if $t \geq t_1 = \ln^2 n + 2\ln n$. In particular, by (15),

$$\mathbf{Pr}\,(M_{t_1} \leq k-1 | X_0 = x) \leq e^{-np_k} + e^{-\ln^2 n}.$$

Since $np_k = o(1)$,

$$e^{-np_k} + e^{-\ln^2 n} \leq e^{-np_k}(1 + 2e^{-\ln^2 n})$$

for $n$ sufficiently large, which we now assume. Thus, for $i = 0, 1, \ldots,$

$$\mathbf{Pr}\,(M_{(i+1)t_1} \leq k-1 | \phi_{it_1}) \leq e^{-np_k}(1 + 2e^{-\ln^2 n})$$

on the event $D_i = (|X_{it_1}| \leq n^2/2) \wedge (M_{it_1} \leq k-1)$. Hence, if we denote $\lfloor n^K/t_1 \rfloor$ by $i_0$, we have

$$\mathbf{Pr}\left((Z \leq k-1) \wedge C\right) \leq \mathbf{Pr}\left(\bigwedge_{i=0}^{i_0} D_i\right)$$

$$= \mathbf{Pr}(D_0) \prod_{i=0}^{i_0-1} \mathbf{Pr}\left(D_{i+1} \Big| \bigwedge_{j=0}^{i} D_j\right)$$

$$\leq (e^{-np_k}(1 + 2e^{-\ln^2 n}))^{i_0}$$

$$\leq (1 + o(1)) \cdot \exp(-(n^K/t_1 - 1)n^{-(K-\varepsilon+o(1))})$$

$$= \exp(-n^{\varepsilon+o(1)}) \to 0$$

as $n \to \infty$. Above we used the observation that

$$(1 + 2e^{-\ln^2 n})^{i_0} \leq \exp(i_0 \cdot 2e^{-\ln^2 n}) = 1 + o(1).$$

**5. Concentration.** We have seen that our balls-and-bins model exhibits rapid mixing. In many Markov models rapid mixing goes along with tight concentration of measure. This is indeed the case here, as demonstrated by the following lemma, which is crucial to our analysis. See [5] for large deviations bounds for a related discrete-time balls-and-bins model.

Let $n$ be a positive integer, and let $\Omega$ be the corresponding set of load vectors, that is, the set of nonnegative vectors in $\mathbb{Z}^n$. A real-valued function $f$ on $\Omega$ is called *Lipschitz* (with Lipschitz constant 1) if

$$|f(x) - f(y)| \leq \|x - y\|_1.$$

LEMMA 5.1. *There is a constant $n_0$ such that, for all $n \geq n_0$, the $n$-bin system has the following property. Let the load vector $Y$ have the equilibrium distribution, and let $f$ be a Lipschitz function on $\Omega$. Then, for each $u \geq n^{1/2} \ln^{3/2} n$,*

$$\mathbf{Pr}\left(|f(Y) - \mathbf{E}[f(Y)]| \geq u\right) \leq e^{-(u^2/n)^{1/3}}.$$

As stated in the Introduction, our primary interest is in the maximum load of a bin. We may deduce easily from the last lemma the following result which we shall use several times.

LEMMA 5.2. *Consider the $n$-bin system in equilibrium. For each positive integer $i$, let $L(i)$ be the random number of bins with at least $i$ balls, at say time $t = 0$, and let $l(i) = \mathbf{E}[L(i)]$. Then*

$$\sup_i \mathbf{Pr}\left(|L(i) - l(i)| \geq n^{1/2} \ln^{3/2} n\right) = O(n^{-1});$$

*for any constant $c > 0$,*

$$\mathbf{Pr}\left(\sup_i |L(i) - l(i)| \geq cn^{1/2} \ln^3 n\right) = e^{-\Omega(\ln^2 n)};$$

*and for each integer $r \geq 2$,*

$$\sup_i \{|\mathbf{E}[L(i)^r] - l(i)^r|\} = O(n^{r-1} \ln^3 n).$$

PROOF. Note that

$$\mathbf{Pr}\left(L(\lceil 2\lambda n \rceil) > 0\right) \leq \mathbf{Pr}\left(\text{Po}(\lambda n) \geq 2\lambda n\right) = e^{-\Omega(n)},$$

since the total number of balls is $\text{Po}(\lambda n)$. Since always $L(i) \leq n$, this shows that we may restrict attention to $i < 2\lambda n$. The first two parts of the lemma now follow directly from Lemma 5.1 (note that $n_0$ is a constant, and does not depend on $f$). To prove the third part, set $u = (r+1)^{3/2} n^{1/2} \ln^{3/2} n$, and note that, by Lemma 5.1,

$$\mathbf{Pr}\left(|L(i) - l(i)| > u\right) \leq e^{-(r+1)\ln n} = n^{-(r+1)}$$

for $n \geq n_0$. Hence, for each positive integer $k \leq r$,

$$\mathbf{E}[|L(i) - l(i)|^k] \leq u^k + n^k \mathbf{Pr}\left(|L(i) - l(i)| > u\right) \leq u^k + o(1).$$

The result now follows from

$$
\begin{aligned}
0 \leq \mathbf{E}[L(i)^r] - l(i)^r &= \sum_{k=2}^{r} \binom{r}{k} \mathbf{E}[(L(i) - l(i))^k] l(i)^{r-k} \\
&\leq \sum_{k=2}^{r} \binom{r}{k} \mathbf{E}[|L(i) - l(i)|^k] n^{r-k} \\
&= O(n^{r-1} \ln^3 n). \qquad \square
\end{aligned}
$$

The next lemma extends the second part of the last lemma, and shows that in equilibrium the number $L_t(i)$ of bins with load at least $i$ at time $t$ is unlikely to move far from its mean value $l(i)$. We show that all the values $L_t(i)$ are likely to stay close to $l(i)$ throughout a polynomial length time interval $[0, \tau]$.

LEMMA 5.3. *Let $K > 0$ be an arbitrary constant, and let $\tau = n^K$. Let $X_0$ be in equilibrium. Then*

$$\mathbf{Pr}\left[\sup_{t \in [0,\tau]} \sup_i |L_t(i) - l(i)| \geq n^{1/2} \ln^3 n\right] = e^{-\Omega(\ln^2 n)}.$$

PROOF. By Lemma 5.2, there exists $\gamma > 0$ such that for all $n$ sufficiently large, for each time $t \geq 0$,

$$\mathbf{Pr}\left(\sup_i |L_t(i) - l(i)| \geq n^{1/2} \ln^3 n/2\right) \leq e^{-\gamma \ln^2 n}.$$

We now let $s = n^{-1/2}$ and $b = 2\lambda n^{1/2}$, and use Lemma 2.1(a), inequality (9). □

The rest of this section is devoted to proving Lemma 5.1. The plan of the proof is as follows. Consider a loads process $(X_t)$, where $X_0 = x_0$ for a suitable load vector $x_0$. (We are most interested in the case $x_0 = \mathbf{0}$.) We shall prove concentration for $X_t$, and later deduce concentration for the equilibrium load vector $Y$.

Note first that the equilibrium load of a bin is stochastically at most $\mathrm{Po}(\lambda d)$. For we can couple the load of a single bin with a process where the arrival rate is always exactly $\lambda d$ and the death rate exactly 1, so that the number of balls in the former is no more than in the latter at all times; and for the latter process, the equilibrium number of balls is $\mathrm{Po}(\lambda d)$.

It will be convenient to limit the maximum load of a bin. Let $b = b(n)$ be an integer at least, say, $4 \ln n / \ln \ln n$—we shall specify a value for $b$ later. Assume that $\max_j x_0(j) \leq b/3$. Let $A_t$ be the event that $M_s \leq b$ for all $0 \leq s \leq t$. If temporarily $\tilde{M}_s$ denotes the maximum load of a process in equilibrium, then, by the time "monotonicity" part of Lemma 3.1, we have

$$\mathbf{Pr}(\overline{A}_t) \leq \mathbf{Pr}(\tilde{M}_s \geq 2b/3 \text{ for some } s \in [0, \tau]).$$

Hence, by (9) in Lemma 2.1(b) and by (6),

$$\mathbf{Pr}(\overline{A}_t) \leq (t+1)(2n)\mathbf{Pr}(\mathrm{Po}(\lambda d) \geq b/3)$$
$$= \exp(\ln(t+1) + \ln n - \tfrac{1}{3}b \ln b + O(b)).$$

It follows that, for $n$ sufficiently large, for each time $t \leq e^b$, say,

(18) $$\mathbf{Pr}(\overline{A}_t) \leq e^{-b \ln b / 13}.$$

In fact, we shall ultimately specify values for $t$ and $b$ so that $t = O(b \ln b)$.

Since loads are rarely large, we can approximate the loads process $(X_t)$ by using only a few of the death processes $F_{j,k}$, namely, those with $k \leq b$, which we call the "low" death processes. In fact, we shall model both the original process and the approximating process, by replacing these low death processes by a combined low death Poisson process with rate $nb$, and a "reaper" process (we omit the "grim"), which at each "toll" of the rate $nb$ Poisson process selects uniformly at random a pair $(j, k)$ where $j \in \{1, \ldots, n\}$ and $k \in \{1, \ldots, b\}$, and behaves as if the process $F_{j,k}$ had "rung." Let $\hat{X}_t$ be the approximating process, which uses only the low death processes. Observe that on $A_t$ we have $\hat{X}_t = X_t$. Since $\mathbf{Pr}\,(\overline{A_t})$ is so small, it will suffice for us to prove concentration for $\hat{X}_t$.

Let $z$ and $\tilde{z}$ be positive integers. Let $\mathbf{t} = (t_1, \ldots, t_z)$ be $z$ arrival times (not ordered) and let $\mathbf{d} = (d_1, \ldots, d_z)$ be corresponding choices of $d$ bins. Let $\tilde{\mathbf{t}} = (\tilde{t}_1, \ldots, \tilde{t}_{\tilde{z}})$ be $\tilde{z}$ low death times (not ordered) and let $\tilde{\mathbf{d}} = (\tilde{d}_1, \ldots, \tilde{d}_{\tilde{z}})$ be corresponding reaper choices [of pairs $(j,k)$, where $j \in \{1, \ldots, n\}$ and $k \in \{1, \ldots, b\}$]. Assume that all these times are distinct. Given any initial load vector $x$, our approximating simulation generates a load vector $s_t(x, \mathbf{t}, \mathbf{d}, \tilde{\mathbf{t}}, \tilde{\mathbf{d}})$ for each time $t \geq 0$.

The following deterministic lemma is analogous to the first part of Lemma 3.1, when the arrivals, choices, death times and reaper choices processes are all deterministic, and may be proved in a similar way.

LEMMA 5.4. *Suppose that we are given two initial load vectors $x_0$ and $y_0$, together with any sequence of arrival times $\mathbf{t}$ and corresponding bin choices $\mathbf{d}$, and departure times $\tilde{\mathbf{t}}$ and corresponding reaper choices $\tilde{\mathbf{d}}$, where all these times are distinct. Then the distance $\|s_t(x_0, \mathbf{t}, \mathbf{d}, \tilde{\mathbf{t}}, \tilde{\mathbf{d}}) - s_t(y_0, \mathbf{t}, \mathbf{d}, \tilde{\mathbf{t}}, \tilde{\mathbf{d}})\|_1$ is nonincreasing in $t$, and so, in particular, for each $t \geq 0$,*

$$\|s_t(x_0, \mathbf{t}, \mathbf{d}, \tilde{\mathbf{t}}, \tilde{\mathbf{d}}) - s_t(y_0, \mathbf{t}, \mathbf{d}, \tilde{\mathbf{t}}, \tilde{\mathbf{d}})\|_1 \leq \|x_0 - y_0\|_1.$$

*Similarly, $\|s_t(x_0, \mathbf{t}, \mathbf{d}, \tilde{\mathbf{t}}, \tilde{\mathbf{d}}) - s_t(y_0, \mathbf{t}, \mathbf{d}, \tilde{\mathbf{t}}, \tilde{\mathbf{d}})\|_\infty$ is nonincreasing in $t$ [recall that $\|z\|_\infty = \max_j |z(j)|$].*

Let us now sketch the plan of the rest of the proof. Let $\mu(t) = \mathbf{E}[f(X_t)]$ and $\hat{\mu}(t) = \mathbf{E}[f(\hat{X}_t)]$. Let $Z_t$ be the number of arrivals in $[0,t]$, so that $Z_t \sim \mathrm{Po}(\lambda n t)$. Let $\tilde{Z}_t$ be the number of low death times in $[0,t]$, so that $\tilde{Z}_t \sim \mathrm{Po}(bnt)$. We shall condition on $Z_t = z$ and $\tilde{Z}_t = \tilde{z}$. Let $\hat{\mu}(t, z, \tilde{z}) = \mathbf{E}[f(\hat{X}_t) | Z_t = z, \tilde{Z}_t = \tilde{z}]$. We shall use Lemma 5.4 and the bounded differences method to upper bound $\mathbf{Pr}\,(|f(\hat{X}_t) - \hat{\mu}(t, z, \tilde{z})| \geq u | Z_t = z, \tilde{Z}_t = \tilde{z})$, see (20) below.

Next we remove the conditioning on $Z_t$ and $\tilde{Z}_t$. To do this, we choose suitable "widths" $w$ and $\tilde{w}$, then use the fact that both $\mathbf{Pr}\,(|Z_t - \lambda n t| > w)$ and $\mathbf{Pr}\,(|\tilde{Z}_t - bnt| > \tilde{w})$ are small, and for $z$ and $\tilde{z}$ such that $|z - \lambda n t| \leq w$ and $|\tilde{z} - bnt| \leq \tilde{w}$, the difference $|\hat{\mu}(t, z, \tilde{z}) - \hat{\mu}(t)|$ is at most about $2(w + \tilde{w})$, see (23) below. We thus find that $\mathbf{Pr}\,(|f(\hat{X}_t) - \hat{\mu}(t)| \geq 3(w + \tilde{w}))$ is small. But since $\hat{X}_t = X_t$ on $A_t$, and $A_t$ is very likely to occur, this last result yields concentration for $f(X_t)$ around its mean. The part of the proof up to here is contained in Lemma 5.5 below. Finally, we use the coupling lemma (Lemma 3.1) to relate the distribution of $X_t$ (with $X_0 = \mathbf{0}$) to the equilibrium distribution.

Let us start on the details. We shall use the following lemma with $x_0 = \mathbf{0}$. (It is convenient elsewhere to have the more general form.)

LEMMA 5.5. *There are constants $n_0$ and $c > 0$ such that the following holds. Let $n \geq n_0$ and $b \geq 4\ln n / \ln \ln n$ be integers, and let $f$ be a Lipschitz function on $\Omega$. Let also $x_0 \in \Omega$ be such that $\max_j x_0(j) \leq b/3$, and assume that the process $(X_t)$ satisfies $X_0 = x_0$ a.s. Then for all times $0 < t \leq e^b$ and all $u \geq 1$,*

$$(19) \qquad \mathbf{Pr}\left(|f(X_t) - \mu_t| \geq u\right) \leq n e^{-cu^2/(nbt)} + e^{-cnt} + e^{-cb\ln b}.$$

PROOF. Note first that we may assume without loss of generality that $f(x_0) = 0$, and so $|f(X_t)| \leq Z_t + \tilde{Z}_t$, since we could replace $f(x)$ by $f(x) - f(x_0)$. Let $z, \tilde{z}$ be positive integers, and condition on $Z_t = z, \tilde{Z}_t = \tilde{z}$. Then $\hat{X}_t$ depends on $2(z + \tilde{z})$ independent random variables $T_1, \ldots, T_z, D_1, \ldots, D_z, \tilde{T}_1, \ldots, \tilde{T}_{\tilde{z}}$, and $\tilde{D}_1, \ldots, \tilde{D}_{\tilde{z}}$. Indeed, we may write $\hat{X}_t$ as $s_t(x_0, \mathbf{T}, \mathbf{D}, \tilde{\mathbf{T}}, \tilde{\mathbf{D}})$, where $\mathbf{T} = (T_1, \ldots, T_z)$, $\mathbf{D} = (D_1, \ldots, D_z)$, $\tilde{\mathbf{T}} = (\tilde{T}_1, \ldots, \tilde{T}_{\tilde{z}})$, and $\tilde{\mathbf{D}} = (\tilde{D}_1, \ldots, \tilde{D}_{\tilde{z}})$. This property relies on the well-known fact that, conditional on the number of events of a Poisson process during $[0, t]$, the unordered event times are a sample of i.i.d. random variables uniform on $[0, t]$. Write

$$g(\mathbf{t}, \mathbf{d}, \tilde{\mathbf{t}}, \tilde{\mathbf{d}}) = f(s_t(x_0, \mathbf{t}, \mathbf{d}, \tilde{\mathbf{t}}, \tilde{\mathbf{d}})).$$

We prove that, conditional on $Z_t = z$ and $\tilde{Z}_t = \tilde{z}$, the random variable $f(\hat{X}_t)$ is highly concentrated, by showing that $g$ satisfies a "bounded differences" condition. Suppose first that we alter a single co-ordinate value $d_j$. Then the value of $g$ can change by at most 2, by Lemma 5.4 starting at time $t_j$ with $\|x_{t_j} - y_{t_j}\|_1 \leq 2$; the same holds if we alter a single co-ordinate value $\tilde{d}_j$. Similarly, if we change a co-ordinate value $t_j$ or $\tilde{t}_j$, the value of $g$ can change by at most 2: we may see this by applying Lemma 5.4 once at the earlier time and once at the later time. Thus, changing any one of the $2(z + \tilde{z})$ co-ordinates can change the value of $g$ by at most 2. Now we use the independent bounded differences inequality, see, for instance, [14]. We find that, for each $u > 0$,

$$\mathbf{Pr}\left(|g(\mathbf{T}, \mathbf{D}, \tilde{\mathbf{T}}, \tilde{\mathbf{D}}) - \mathbf{E}[g(\mathbf{T}, \mathbf{D}, \tilde{\mathbf{T}}, \tilde{\mathbf{D}})]| \geq u\right) \leq 2\exp\left(-\frac{u^2}{4(z + \tilde{z})}\right).$$

In other words, we have proved that, for any $u > 0$,

$$(20) \quad \mathbf{Pr}\left(|f(\hat{X}_t) - \hat{\mu}(t, z, \tilde{z})| \geq u | Z_t = z, \tilde{Z}_t = \tilde{z}\right) \leq 2\exp\left(-\frac{u^2}{4(z + \tilde{z})}\right).$$

Next we will remove the conditioning on $Z_t$ and $\tilde{Z}_t$. We will choose suitable "widths" $w = w(n)$ and $\tilde{w} = \tilde{w}(n)$, where $0 \leq w \leq \lambda nt$ and $0 \leq \tilde{w} \leq bnt$. Let $I$ denote the interval of integer values $z$ such that $|z - \lambda nt| \leq w$, and let $\tilde{I}$ denote the interval of integer values $\tilde{z}$ such that $|\tilde{z} - bnt| \leq \tilde{w}$. Recall that we shall ensure that with high probability $Z_t \in I$ and $\tilde{Z}_t \in \tilde{I}$, and for each $z \in I$ and $\tilde{z} \in \tilde{I}$, the difference $|\hat{\mu}(t, z, \tilde{z}) - \hat{\mu}(t)|$ is not too large.

Since $Z_t \sim \text{Po}(\lambda nt)$ and $\tilde{Z}_t \sim \text{Po}(bnt)$, by (5),

$$\textbf{Pr}\,(Z_t \notin I) = \textbf{Pr}\,(|Z_t - \lambda nt| > w) \leq 2\exp\left(-\frac{w^2}{3\lambda nt}\right) \tag{21}$$

and

$$\textbf{Pr}\,(\tilde{Z}_t \notin \tilde{I}) = \textbf{Pr}\,(|\tilde{Z}_t - bnt| > \tilde{w}) \leq 2\exp\left(-\frac{\tilde{w}^2}{3bnt}\right). \tag{22}$$

We shall choose $w$ and $\tilde{w}$ to satisfy $w \geq 2(\lambda nt \ln n)^{1/2}$ and $\tilde{w} \geq 2(bnt \ln n)^{1/2}$. Then, by (21), (22), (5) and (7), provided that $b$ satisfies $b = o(n^{1/3})$,

$$\textbf{E}[Z_t \mathbf{1}_{(Z_t \notin I \vee \tilde{Z}_t \notin \tilde{I})}] \leq \textbf{E}[Z_t \mathbf{1}_{Z_t > \lambda nt + w}] + \lambda nt \,\textbf{Pr}\,(\tilde{Z}_t \notin \tilde{I}) = o(1)$$

and

$$\textbf{E}[\tilde{Z}_t \mathbf{1}_{(Z_t \notin I \vee \tilde{Z}_t \notin \tilde{I})}] \leq \textbf{E}[\tilde{Z}_t \mathbf{1}_{\tilde{Z}_t > bnt + \tilde{w}}] + bnt \,\textbf{Pr}\,(Z_t \notin I) = o(1).$$

Hence, since $|f(\hat{X}_t)| \leq Z_t + \tilde{Z}_t$,

$$\textbf{E}[|f(\hat{X}_t)| \, \mathbf{1}_{(Z_t \notin I \vee \tilde{Z}_t \notin \tilde{I})}] = o(1).$$

But

$$\hat{\mu}(t) = \sum_{z \in I, \tilde{z} \in \tilde{I}} \hat{\mu}(t, z, \tilde{z}) \,\textbf{Pr}\,(Z_t = z, \tilde{Z}_t = \tilde{z}) + \textbf{E}[f(\hat{X}_t) \mathbf{1}_{(Z_t \notin I \vee \tilde{Z}_t \notin \tilde{I})}].$$

Hence,

$$\hat{\mu}(t) \leq \max_{z \in I, \tilde{z} \in \tilde{I}} \{\hat{\mu}(t, z, \tilde{z})\} + o(1),$$

and, using also (21) and (22),

$$\hat{\mu}(t) \geq \min_{z \in I, \tilde{z} \in \tilde{I}} \{\hat{\mu}(t, z, \tilde{z})\} + o(1).$$

By Lemma 5.4, for each $z, \tilde{z}$,

$$|\hat{\mu}(t, z+1, \tilde{z}) - \hat{\mu}(t, z, \tilde{z})| \leq 1$$

and

$$|\hat{\mu}(t, z, \tilde{z}+1) - \hat{\mu}(t, z, \tilde{z})| \leq 1.$$

It follows, using the bounds above on $\hat{\mu}(t)$, that, for each $z \in I$ and $\tilde{z} \in \tilde{I}$,

$$|\hat{\mu}(t, z, \tilde{z}) - \hat{\mu}(t)| \leq 2(w + \tilde{w}) + o(1). \tag{23}$$

Now by (20), (21), (22) and (23),

$$\mathbf{Pr}\left(|f(\hat{X}_t) - \hat{\mu}(t)| \geq (3 + o(1))(w + \tilde{w})\right)$$
$$\leq \sum_{z \in I, \tilde{z} \in \tilde{I}} \mathbf{Pr}\left(|f(\hat{X}_t) - \hat{\mu}(t)| \geq (3 + o(1))(w + \tilde{w}) | Z_t = z, \tilde{Z}_t = \tilde{z}\right)$$
$$\times \mathbf{Pr}\left(Z_t = z, \tilde{Z}_t = \tilde{z}\right)$$
$$+ \mathbf{Pr}\left(Z_t \notin I\right) + \mathbf{Pr}\left(\tilde{Z}_t \notin \tilde{I}\right)$$
$$\leq \sum_{z \in I, \tilde{z} \in \tilde{I}} \mathbf{Pr}\left(|f(\hat{X}_t) - \hat{\mu}(t, z, \tilde{z})| \geq (1 + o(1))(w + \tilde{w}) | Z_t = z, \tilde{Z}_t = \tilde{z}\right)$$
$$\times \mathbf{Pr}\left(Z_t = z, \tilde{Z}_t = \tilde{z}\right)$$
$$+ \mathbf{Pr}\left(|Z_t - \lambda nt| > w\right) + \mathbf{Pr}\left(|\tilde{Z}_t - bnt| > \tilde{w}\right)$$
$$\leq 2 \exp\left(-\frac{(1 + o(1))(w + \tilde{w})^2}{4(\lambda nt + bnt + w + \tilde{w})}\right) + 2 \exp\left(-\frac{w^2}{3\lambda nt}\right) + 2 \exp\left(-\frac{\tilde{w}^2}{3bnt}\right)$$
$$\leq 2 \exp\left(-\frac{(1 + o(1))(w + \tilde{w})^2}{5nbt}\right) + 2 \exp\left(-\frac{w^2}{3\lambda nt}\right) + 2 \exp\left(-\frac{\tilde{w}^2}{3bnt}\right),$$

since $b(n) \to \infty$ as $n \to \infty$. Let $u$ satisfy

$$6(nbt \ln n)^{1/2} \leq u \leq 3\sqrt{\lambda b}nt.$$

Let $\tilde{w} = u/3$ and $w = \tilde{w}\sqrt{\lambda/b}$. Observe that, for $n$ sufficiently large, the bounds required above on $w$ and $\tilde{w}$ hold, and $u = (3 + o(1))(w + \tilde{w})$. Thus,

$$\mathbf{Pr}\left(|f(\hat{X}_t) - \hat{\mu}(t)| \geq u\right) \leq 2e^{-(1+o(1))u^2/(45nbt)} + 4e^{-u^2/(27nbt)}$$
$$\leq e^{-u^2/(46nbt)}$$

for $n$ sufficiently large. But if $u < 6(nbt \ln n)^{1/2}$, then $e^{-u^2/(46nbt)} \geq n^{-1}$. Thus, as long as $u \leq 3\sqrt{\lambda b}nt$, we have

$$\mathbf{Pr}\left(|f(\hat{X}_t) - \hat{\mu}(t)| \geq u\right) \leq n e^{-u^2/(46nbt)}.$$

Now we move from $\hat{X}_t$ to $X_t$. Note that in $[0, t]$ there are $Z_t$ arrivals and at most $|X_0| + Z_t$ departures, and so $|f(X_t) - f(\hat{X}_t)| \leq 2(|X_0| + 2Z_t)$. Thus, since also $X_t = \hat{X}_t$ on $A_t$,

$$|\hat{\mu}(t) - \mu(t)| = |\mathbf{E}[(f(X_t) - f(\hat{X}_t))\mathbf{1}_{\overline{A}_t}]| \leq 2\mathbf{E}[(|X_0| + 2Z_t)\mathbf{1}_{\overline{A}_t}].$$

But $|X_0| \leq nb/3$ and $\mathbf{E}[Z_t \mathbf{1}_{\overline{A}_t}] \leq 2\lambda nt \mathbf{Pr}\left(\overline{A}_t\right) + \mathbf{E}[Z_t \mathbf{1}_{Z_t > 2\lambda nt}]$. Hence,

$$|\hat{\mu}(t) - \mu(t)| \leq (2nb/3 + 8\lambda nt)\mathbf{Pr}\left(\overline{A}_t\right) + 4\mathbf{E}[Z_t \mathbf{1}_{Z_t > 2\lambda nt}] = o(1),$$

by (18) and (7). Thus,

$$\mathbf{Pr}\left(|f(X_t) - \mu(t)| \geq u\right) \leq \mathbf{Pr}\left(|f(\hat{X}_t) - \hat{\mu}(t)| \geq u + o(1)\right) + \mathbf{Pr}\left(\overline{A_t}\right)$$
$$\leq ne^{-(u+o(1))^2/(46nbt)} + e^{-b\ln b/13},$$

by (18) (since we assume that $t \leq e^b$). The lemma now follows, by replacing $u$ by $\min\{u, 3\sqrt{\lambda bnt}\}$. $\square$

We shall use Lemma 5.5 here with $X_0 = \mathbf{0}$ to complete the proof of Lemma 5.1. As we saw before, we may assume that $f(\mathbf{0}) = 0$, and, hence, always $|f(x)| \leq |x|$. It remains to relate the distribution of $X_t$ with $X_0 = \mathbf{0}$ to the equilibrium distribution, and to choose values for $b$ and $t$. By Theorem 1.1, if $Y$ has the equilibrium distribution, then $d_{\mathrm{TV}}(\mathcal{L}(X_t), \mathcal{L}(Y)) \leq \lambda n e^{-t}$. Hence, for all $n$ sufficiently large, $b \geq 4 \ln n / \ln \ln n$ and $u \geq 1$,

$$\mathbf{Pr}\left(|f(Y) - \mu(t)| \geq u\right)$$
(24)
$$\leq d_{\mathrm{TV}}(\mathcal{L}(X_t), \mathcal{L}(Y)) + \mathbf{Pr}\left(|f(X_t) - \mu(t)| \geq u\right)$$
$$\leq \lambda n e^{-t} + n e^{-cu^2/(ntb)} + e^{-cnt} + e^{-cb\ln b}.$$

Let $u \geq 2(n \ln^3 n / c \ln \ln n)^{1/2}$. Let $t = (u^2 c \ln \ln n / n)^{1/3}$ and $b = \lfloor t / \ln \ln n \rfloor$. Then $t \geq 4^{1/3} \ln n$. Also, $\ln b \geq (1 + o(1)) \ln \ln n$, so $b \ln b \geq (1 + o(1))t = \Omega(t)$. Further, $cu^2/(nbt) = \Omega(t)$. It now follows from (24) that

$$\mathbf{Pr}\left(|f(Y) - \mu(t)| \geq u\right) = e^{-\Omega((u^2 \ln \ln n / n)^{1/3})}.$$

Finally, we relate $\mu(t) = \mathbf{E}[f(X_t)]$ to $\mathbf{E}[f(Y)]$. By Theorem 1.1,

$$|\mu(t) - \mathbf{E}[f(Y)]| \leq d_{\mathrm{W}}(\mathcal{L}(X_t), \mathcal{L}(Y)) = o(1)$$

since $t \geq 4^{1/3} \ln n$. Thus, we find that, for any $u \geq 2(n \ln^3 n / c \ln \ln n)^{1/2}$,

$$\mathbf{Pr}\left(|f(Y) - \mathbf{E}[f(Y)]| \geq u\right) = e^{-\Omega((u^2 \ln \ln n / n)^{1/3})}.$$

This completes the proof of Lemma 5.1.

**6. Balance equations.** In this section we suppose throughout that the system is in equilibrium. We present the balance equation (26), and deduce Lemma 6.1, which we shall need in Section 8, concerning the expected proportion of bins with at least $i$ balls.

Let $d \geq 2$ be a fixed integer. Consider a positive integer $n$, and the corresponding set $\Omega$ of load vectors. For $x \in \Omega$ and a nonnegative integer $k$, let $u(k, x)$ be the proportion of bins $j$ with load $x(j)$ at least $k$. Thus, always $u(0, x) = 1$. Let $X$ have the equilibrium distribution over $\Omega$, and let $u(k)$ denote $\mathbf{E}[u(k, X)]$ (which depends on $n$). [Thus, $u(k) = l(k)/n$, where $l(k)$ was defined earlier as the expected number of bins with load at least $k$.]

LEMMA 6.1.  (a) *There is a constant $c$ such that, for $n$ sufficiently large, if $j \geq \ln\ln n / \ln d + c$, then $u(j) \leq n^{-1} \ln^3 n$.*

(b) *For any $\eta > 0$, there is a constant $c$ such that, for $n$ sufficiently large, if $j \leq \ln\ln n / \ln d - c$, then $u(j) \geq n^{-\eta}$.*

The rest of this section is devoted to proving this lemma. First we present the balance equations.

It is easy to check (see [21]) that, if $f$ is the bounded function $f(x) = u(k,x)$, then the generator operator $G$ of the Markov process satisfies

$$Gf(x) = \lambda(u(k-1,x)^d - u(k,x)^d) - k(u(k,x) - u(k+1,x))$$

[cf. with equation (1) earlier]. To see this, note that $u(k,x) - u(k+1,x)$ is the proportion of bins with load exactly $k$, and $u(k-1,x)^d - u(k,x)^d$ is the probability that the minimum load of the $d$ attempts is exactly $k-1$. Since $X$ is in equilibrium, $\mathbf{E}[Gf(X)] = 0$. Hence,

(25)  $\lambda(\mathbf{E}[u(k-1,X)^d] - \mathbf{E}[u(k,X)^d]) - k(u(k) - u(k+1)) = 0.$

Now

$$\sum_{k \geq 1} k u(k,x) = \frac{1}{n} \sum_{j=1}^{n} \binom{x(j)+1}{2} \leq \frac{|x|^2}{n},$$

and so

$$\sum_{k \geq 1} k u(k) \leq \frac{\mathbf{E}[|X|^2]}{n} < \infty.$$

Hence, $ku(k) \to 0$ as $k \to \infty$. Also, $\mathbf{E}[u(k,X)^d] \leq u(k)$. It follows on summing (25), for $k \geq i$, that, for each $i = 1, 2 \ldots$, we have

(26)  $\lambda \mathbf{E}[u(i-1,X)^d] - i u(i) - \sum_{k \geq i+1} u(k) = 0.$

(This is the result that $\mathbf{E}[Gf(X)] = 0$, where $f(x)$ is the number of balls of "height" at least $i$, i.e., $f(x) = \sum_{j=1}^{n}(x(j) - i + 1)^+$, but since $f$ is not bounded, we cannot assert the result directly.) Equation (26) is the key fact in our analysis. Observe that, by (26), for each positive integer $i$,

(27)  $$u(i) \leq \frac{\lambda}{i} \mathbf{E}[u(i-1,X)^d].$$

We are now ready to prove the lemma, part (b) first. Let $a = \lceil 2\lambda \rceil - 1$. We shall show that $u(a)$ is at least a positive constant, and the $u(i)$ do not decrease too quickly for $i \geq a$.

Note first that, since $\mathbf{E}[u(i-1,X)^d] \geq u(i-1)^d$, by (26), we have

(28)  $$\lambda u(i-1)^d - i u(i) - \sum_{k \geq i+1} u(k) \leq 0.$$

Also, since $0 \leq u(i-1, X) \leq 1$, we have $\mathbf{E}[u(i-1, X)^d] \leq u(i-1)$ and so by (27), for each $i = 1, 2, \ldots$, we have $u(i) \leq \lambda u(i-1)/i$. Thus, for $i \geq a$, we have $u(i+1) \leq \lambda u(i)/(i+1) \leq u(i)/2$. Hence, if $k \geq i \geq a$, then $u(k) \leq 2^{-(k-i)} u(i)$; and so

$$\sum_{k \geq i+1} u(k) \leq u(i) \qquad \text{for } i \geq a.$$

It now follows from (28) that, for $i \geq a$, we have $\lambda u(i-1)^d - (i+1)u(i) \leq 0$; and, thus,

$$(29) \qquad u(i) \geq \frac{\lambda u(i-1)^d}{i+1} \qquad \text{for } i \geq a.$$

Inequality (29) will show that the $u(i)$ do not decrease too quickly for $i \geq a$.

Now consider small values of $i$. Let $i \in \{1, \ldots, a\}$. Since $u(i) \geq u(k)$ for $k \geq i$, we have $(a-i)u(i) - \sum_{k=i+1}^{a} u(k) \geq 0$. Hence, by (28),

$$0 \geq \lambda u(i-1)^d - iu(i) - \sum_{k \geq i+1} u(k)$$
$$\geq \lambda u(i-1)^d - au(i) - \sum_{k \geq a+1} u(k)$$
$$\geq \lambda u(i-1)^d - (a+1)u(i).$$

Thus, we have

$$u(i) \geq \frac{\lambda}{2\lambda+1} u(i-1)^d \qquad \text{for } i = 1, \ldots, a.$$

The last inequality shows that there is a constant $\delta_1 > 0$ (depending only on $\lambda$ and $d$) such that always $u(a) \geq \delta_1$. But by (29) and induction on $i$, for each $i = 1, 2, \ldots$,

$$u(a+i) \geq \frac{\lambda^{1+d+\cdots+d^{i-1}}}{(a+i+1)(a+i)^d(a+i-1)^{d^2} \cdots (a+2)^{d^{i-1}}} u(a)^{d^i}.$$

To upper bound the denominator, note that

$$\ln((a+i+1)(a+i)^d(a+i-1)^{d^2} \cdots (a+2)^{d^{i-1}})$$
$$= d^i \sum_{k=1}^{i} d^{-k} \ln(a+k+1) \leq c_2 \, d^i$$

for some constant $c_2$, and so the denominator is at most $e^{c_2 d^i}$. Let $\delta_2 > 0$ be the constant $\lambda e^{-c_2} \delta_1$. Then

$$u(i) \geq u(a+i) \geq \delta_2^{d^i} = \exp\left(-d^i \ln \frac{1}{\delta_2}\right)$$

for each $i = 1, 2, \ldots$. Let the constant $c_3$ be such that $d^{-c_3} \ln \frac{1}{\delta_2} \leq \eta$. Hence, if $i \leq \ln\ln n / \ln d - c_3$, then

$$u(i) \geq \exp\left(-(\ln n)d^{-c_3} \ln \frac{1}{\delta_2}\right)$$

$$\geq \exp(-\eta \ln n) = n^{-\eta}.$$

This completes the proof of part (b) of the lemma.

We now prove part (a) of the lemma. By Lemma 5.2, there exists a constant $c_1 > 0$ such that, for all positive integers $i$ and $n$,

(30) $$u(i) \leq \frac{\lambda}{i}(u(i-1)^d + c_1 n^{-1} \ln^3 n).$$

Let $i^* = i^*(n)$ be the smallest positive integer $i$ such that $u(i-1)^d < c_1 n^{-1} \ln^3 n$, that is, $u(i-1) < c_1^{1/d} n^{-1/d} (\ln n)^{3/d}$. We may assume that $n$ is sufficiently large that $c_1 \ln^3 n > 1$, and so the quantity $c_1^{1/d} n^{-1/d} (\ln n)^{3/d}$ in the last bound is $> 1/n$. Note that, by (30),

$$u(i^*) \leq \frac{2\lambda}{i^*} c_1 n^{-1} \ln^3 n = o(n^{-1} \ln^3 n),$$

since $i^*(n) \to \infty$ as $n \to \infty$ by part (b).

We want an upper bound on $i^*$. By (30),

(31) $$u(i) \leq \frac{2\lambda}{i} u(i-1)^d$$

for each $i = 1, \ldots, i^* - 1$. Let $i_0$ be the constant $\lceil 2e\lambda \rceil$. We check that $i^* < \ln\ln n / \ln d + i_0 + 2$. Since $i^*(n) \to \infty$ as $n \to \infty$, we may assume that $i_0 \leq i^* - 1$. By (31), $u(i_0) \leq \frac{2\lambda}{i_0} u(i_0 - 1)^d \leq \frac{2\lambda}{i_0} \leq e^{-1}$. Also by (31), for $i = i_0 + 1, \ldots, i^* - 1$, we have $u(i) \leq u(i-1)^d$, and it follows that $u(i) \leq e^{-d^{i-i_0}}$ for each $i = i_0, \ldots, i^* - 1$. But $e^{-d^{i-i_0}} \leq 1/n$ when $d^{i-i_0} \geq \ln n$; that is, when $i \geq \ln\ln n / \ln d + i_0$. Thus, if $i^* \geq \ln\ln n / \ln d + i_0 + 2$, then $u(i^* - 2) \leq 1/n$, contradicting the choice of $i^*$. This completes the proof of part (a) of Lemma 6.1, and thus of the whole lemma.

**7. Random walks with drift.** In this section we consider a generalized random walk on the integers, which takes steps of $0, \pm 1$ but with probabilities that can depend on the history of the process, and where there is a "drift." We shall use the following version of the Bernstein inequality—see Theorem 2.7 in [14].

LEMMA 7.1. *Let $b \geq 0$, and let the random variables $Z_1, \ldots, Z_n$ be independent, with $Z_k - \mathbf{E}[Z_k] \geq -b$ for each $k$. Let $S_n = \sum_k Z_k$, and let $S_n$ have* *expected value $\mu$ and variance $V$ (assumed finite). Then for any $z \geq 0$,*

$$\mathbf{Pr}\left(S_n \leq \mu - z\right) \leq \exp\left(-\frac{z^2}{2V + (2/3)bz}\right). \tag{32}$$

(The term $\frac{2}{3}bz$ should be thought of as an error term.) The next lemma concerns hitting times for a generalized random walk with "drift."

LEMMA 7.2. *Let $\phi_0 \subseteq \phi_1 \subseteq \cdots \subseteq \phi_m$ be a filtration, and let $Y_1, Y_2, \ldots, Y_m$ be random variables taking values in $\{-1, 0, 1\}$ such that each $Y_i$ is $\phi_i$-measurable. Let $E_0, E_1, \ldots, E_{m-1}$ be events, where $E_i \in \phi_i$ for each $i$, and let $E = \bigwedge_i E_i$. Let $R_t = R_0 + \sum_{i=1}^{t} Y_i$. Let $0 \leq p \leq 1/3$, let $r_0$ and $r_1$ be integers such that $r_0 < r_1$, and let $pm \geq 2(r_1 - r_0)$. Assume that, for each $i = 1, \ldots, m$,*

$$\mathbf{Pr}\left(Y_i = 1 | \phi_{i-1}\right) \geq 2p \qquad \text{on } E_{i-1} \wedge (R_{i-1} < r_1)$$

*and*

$$\mathbf{Pr}\left(Y_i = -1 | \phi_{i-1}\right) \leq p \qquad \text{on } E_{i-1} \wedge (R_{i-1} < r_1).$$

*Then*

$$\mathbf{Pr}\left(E \wedge (R_t < r_1 \ \forall t \in \{1, \ldots, m\}) | R_0 = r_0\right) \leq \exp\left(-\frac{pm}{28}\right).$$

PROOF. Let us first prove the lemma assuming that the above inequalities on $\mathbf{Pr}(Y_i = 1 | \phi_{i-1})$ and $\mathbf{Pr}(Y_i = -1 | \phi_{i-1})$ hold a.s.; that is, ignoring the events $E_{i-1} \wedge (R_{i-1} < r_1)$. We shall then see easily how to incorporate these events.

We can couple the $Y_i$ with i.i.d. random variables $Z_i$ taking values in $\{-1, 0, 1\}$, such that $\mathbf{Pr}(Z_i = 1) = 2p$, $\mathbf{Pr}(Z_i = -1) = p$ and $\mathbf{Pr}(Z_i \leq Y_i) = 1$ for each $i$. The variables $Z_1, Z_2, \ldots$ are independent; $\mathbf{E}[Z_i] = p$, $\mathbf{Var}[Z_i] \leq 3p$, and $Z_i - \mathbf{E}[Z_i] \geq -1 - p \geq -4/3$ for each $i$. Let $S_t = \sum_{i=1}^{t} Z_i$, let $\mu_t = \mathbf{E}(S_t) = pt$, and note that $\mathbf{Var}(S_t) \leq 3tp$. Hence, by Bernstein's inequality, Lemma 7.1, for each $y > 0$,

$$\mathbf{Pr}\left(S_t \leq \mu_t - y\right) \leq \exp\left(-\frac{y^2}{6pt + y}\right).$$

Note that $\mu_m = pm$. Thus, if $a = r_1 - r_0$,

$$\mathbf{Pr}\left(R_t < r_1 \ \forall t \in \{1, \ldots, m\} | R_0 = r_0\right)$$
$$\leq \mathbf{Pr}\left(S_m < a\right)$$
$$\leq \exp\left(-\frac{(pm - a)^2}{6pm + (pm - a)}\right)$$

$$\leq \exp\left(-\frac{pm}{7}\left(1-\frac{a}{pm}\right)^2\right)$$

$$\leq \exp\left(-\frac{pm}{28}\right),$$

since $a/pm \leq 1/2$.

Now let us return to the full lemma as stated, with the events $E_i$. For each $i = 0, 1, \ldots, m-1$, let $F_i = E_i \wedge (R_i < r_1)$; and for each $i = 1, \ldots, m$, let $\tilde{Y}_i = Y_i \cdot \mathbf{1}_{F_{i-1}} + \mathbf{1}_{\overline{F}_{i-1}}$. Let $\tilde{R}_0 = R_0$ and for $t = 1, \ldots, m$, let $\tilde{R}_t = R_0 + \sum_{i=1}^{t} \tilde{Y}_i$. Then $\mathbf{Pr}(\tilde{Y}_i = 1|\phi_{i-1}) \geq 2p$, since, by assumption, it is at least $2p$ on $F_{i-1}$, and it equals 1 on $\overline{F}_{i-1}$. Similarly, $\mathbf{Pr}(\tilde{Y}_i = -1|\phi_{i-1}) \leq p$. Hence, by what we have just proved applied to the $\tilde{Y}_i$,

$$\mathbf{Pr}(E \wedge (R_t < r_1 \ \forall t \in \{1, \ldots, m\})|R_0 = r_0)$$
$$= \mathbf{Pr}(E \wedge (\tilde{R}_t < r_1 \ \forall t \in \{1, \ldots, m\})|\tilde{R}_0 = r_0)$$
$$\leq \mathbf{Pr}(\tilde{R}_t < r_1 \ \forall t \in \{1, \ldots, m\}|\tilde{R}_0 = r_0)$$
$$\leq \exp\left(-\frac{pm}{28}\right),$$

as required. $\square$

The next lemma shows that if we try to cross an interval against the drift, then we will rarely succeed.

LEMMA 7.3. *Let $a$ be a positive integer. Let $p$ and $q$ be reals with $q > p \geq 0$ and $p + q \leq 1$. Let $\phi_0 \subseteq \phi_1 \subseteq \phi_2 \subseteq \cdots$ be a filtration, and let $Y_1, Y_2, \ldots$ be random variables taking values in $\{-1, 0, 1\}$ such that each $Y_i$ is $\phi_i$-measurable. Let $E_0, E_1, \ldots$ be events where each $E_i \in \phi_i$, and let $E = \bigwedge_i E_i$. Let $R_0 = 0$ and let $R_k = \sum_{i=1}^{k} Y_i$ for $k = 1, 2, \ldots$. Assume that, for each $i = 1, 2, \ldots$,*

$$\mathbf{Pr}(Y_i = 1|\phi_{i-1}) \leq p \qquad \text{on } E_{i-1} \wedge (0 \leq R_{i-1} \leq a-1)$$

*and*

$$\mathbf{Pr}(Y_i = -1|\phi_{i-1}) \geq q \qquad \text{on } E_{i-1} \wedge (0 \leq R_{i-1} \leq a-1).$$

*Let*

$$T = \inf\{k \geq 1 : R_k \in \{-1, a\}\}.$$

*Then*

$$\mathbf{Pr}(E \wedge (R_T = a)) \leq (p/q)^a.$$

24 M. J. LUCZAK AND C. MCDIARMIDPROOF. As with the previous lemma, let us first prove this lemma assuming that the given inequalities on $\mathbf{Pr}(Y_i = 1|\phi_{i-1})$ and $\mathbf{Pr}(Y_i = -1|\phi_{i-1})$ hold a.s. We can couple the $Y_i$ with i.i.d. random variables $\hat{Y}_i$ taking values in $\{0, \pm 1\}$ such that $\mathbf{Pr}(\hat{Y}_i = 1) = p$, $\mathbf{Pr}(\hat{Y}_i = -1) = q$ and $\mathbf{Pr}(Y_i \leq \hat{Y}_i) = 1$. Let $\hat{R}_0 = 0$, let $\hat{R}_k = \sum_{i=1}^{k} \hat{Y}_i$ for $k = 1, 2, \ldots$, and let

$$\hat{T} = \inf\{k \geq 1 : \hat{R}_k \in \{-1, a\}\}.$$

Then from standard properties of a simple random walk,

$$\mathbf{Pr}(R_T = a) \leq \mathbf{Pr}(\hat{R}_{\hat{T}} = a) \leq (p/q)^a.$$

Now let us incorporate the events $E_i$, and consider the full lemma as stated. For each $i = 0, 1, \ldots$, let $F_i = E_i \wedge (0 \leq R_{i-1} \leq a - 1)$; and for each $i = 1, 2, \ldots$, let $\tilde{Y}_i = Y_i \mathbf{1}_{F_{i-1}} - \mathbf{1}_{\overline{F}_{i-1}}$. Let $\tilde{R}_k$ and $\tilde{T}$ be defined in the obvious way. Then $\mathbf{Pr}(\tilde{Y}_i = 1|\phi_{i-1}) \leq p$, since, by assumption, it is at most $p$ on $F_{i-1}$, and it equals 0 on $\overline{F}_{i-1}$. Similarly, $\mathbf{Pr}(\tilde{Y}_i = -1|\phi_{i-1}) \geq q$. Hence, by what we have just proved applied to the $\tilde{Y}_i$,

$$\mathbf{Pr}(E \wedge (R_T = a)) \leq \mathbf{Pr}(\tilde{R}_{\tilde{T}} = a) \leq (p/q)^a. \qquad \square$$

**8. Proof of Theorem 1.3.** We have assembled all the preliminary results we need. In this section we at last prove Theorem 1.3, and inequalities (3) and (4) that follow it. We assume throughout that the process is in equilibrium.

Let $d \geq 2$ be a fixed integer. Consider the $n$-bin system. Recall that $u(k)$ is the expected proportion of bins with load at least $k$. Define $j^* = j^*(n)$ to be the least positive integer $i$ such that $u(i) < n^{-1/2} \ln^3 n$. By Lemma 6.1,

$$j^*(n) = \ln \ln n / \ln d + O(1).$$

We shall show that,

(33) for $d = 2$, we have $M = j^*$ or $j^* + 1$ a.a.s.

and

(34) for $d \geq 3$, we have $M = j^* - 1$ or $j^*$ a.a.s.

This will complete the proof of the first part of Theorem 1.3.

For each time $t$ and each $i = 0, 1, \ldots$, let the random variable $Z_t(i)$ be the number of new balls arriving during $[0, t]$ which have height at least $i$ on arrival, that is, which are placed in a bin already holding at least $i - 1$ balls. Let $J_0 = 0$ and enumerate the arrival times after time 0 as $J_1, J_2, \ldots$. We shall define a "horizon" time $t_0$ of the order of $\ln n$, and let $N = \lceil 2\lambda n t_0 \rceil$. For each time $t$, let $A_t$ be the event

$$\{\lambda n/2 \leq |X_s| \leq 2\lambda n \ \forall \, s \in [0, t]\}.$$

Then by Lemma 2.2, the event $A_{t_0}$ holds with probability $1 - e^{-\Omega(n)}$.

8.1. *The case $d \geq 3$.* We consider first the case when $d \geq 3$, which is easier than when $d = 2$.

Let $K > 0$ be a (large) constant and let $t_0 = (K+4)\ln n$. Since $l(j^* - 1) \geq n^{1/2} \ln^3 n$, the concentration result Lemma 5.3 shows that $\mathbf{Pr}\,(M < j^* - 1) = e^{-\Omega(\ln^2 n)}$. In particular, $M \geq j^* - 1$ a.a.s., which is "half" of (34). Also, (8) in Lemma 2.1 above [with $s = n^{-(K+4)}$, $a = j^* - 3$ and $b = 1$] shows that

$$(35) \qquad \min\{M_t : 0 \leq t \leq n^K\} \geq j^* - 2 \quad \text{a.a.s.}$$

This result establishes a finer form of the lower bound half of (2). In fact, this half of (2) will follow from (3) which we prove later, so (35) is not needed for our proofs.

Next we shall show that $M \leq j^*$ a.a.s., which is the other half of (34). For $k = 0, 1, \ldots$, let $E_k$ be the event that at time $J_k$ there are no more than $2n^{1/2}\ln^3 n$ bins with at least $j^*$ balls. Then $\mathbf{Pr}\,(\overline{E_k}) = e^{-\Omega(\ln^2 n)}$ by Lemma 5.2, since $l(j^*) < n^{1/2}\ln^3 n$. Consider the ball which arrives at time $J_k$: on $E_{k-1}$, it has probability at most $p_1 = (2n^{-1/2}\ln^3 n)^d$ of falling into a bin with at least $j^*$ balls. Note that

$$(36) \qquad \mathbf{Pr}\,(J_{N+1} \leq t_0) \leq \mathbf{Pr}\,[\text{Po}(\lambda n t_0) \geq 2\lambda n t_0] = e^{-\Omega(n \ln n)}.$$

Also, for each positive integer $r$,

$$\mathbf{Pr}\,(B(N, p_1) \geq r) \leq (Np_1)^r = O((n^{-(d/2-1)}(\ln n)^{3d+1})^r).$$

(Here we are using $B$ to denote a binomial random variable.) Hence, for each positive integer $r$, using Lemma 5.3,

$$\mathbf{Pr}\,[Z_{t_0}(j^*+1) \geq r]$$
$$\leq \mathbf{Pr}\,(B(N, p_1) \geq r) + \mathbf{Pr}\,\left(\bigvee_{k=0}^{N-1} \overline{E_k}\right) + \mathbf{Pr}\,(J_{N+1} \leq t_0)$$
$$= O((n^{-(d/2-1)}(\ln n)^{3d+1})^r).$$

Also, the probability that some "initial" ball survives to time $t_0$ is at most $\lambda n e^{-t_0}$, as we saw earlier. Hence, for each positive integer $r$,

$$\mathbf{Pr}\,[M \geq j^* + r] \leq \mathbf{Pr}\,[Z_{t_0}(j^*+1) \geq r] + \lambda n e^{-t_0}.$$

In particular, $\mathbf{Pr}\,(M \geq j^* + 1) = o(1)$, which together with the earlier result that $M \geq j^* - 1$ a.a.s., completes the proof of (34). Further,

$$\mathbf{Pr}\,[M \geq j^* + 2K + 5] = o(n^{-K-2}).$$

Now (9) in Lemma 2.1 with $\tau = n^K$ and $s = n^{-2}$, together with (35), lets us complete the proof of (2).

8.2. *The case $d=2$.* The case $d=2$ needs a little more effort, and uses the "drift" results from the last section. Again, let $K > 0$ be a (large) constant, but now let $t_0 = (2K+8)\ln n$. We first show that $M \geq j^*$, by showing that, in fact,

(37) $$L_{t_0}(j^*) \geq \ln^3 n \quad \text{a.a.s.}$$

Let $J'_0 = 0$, and enumerate all jump times after time 0 (not just the arrival times) as $J'_1, J'_2, \ldots$. Note that $J'_n \leq t_0$ a.a.s., since

$$\mathbf{Pr}(J'_n > t_0) \leq \mathbf{Pr}(J_n > t_0) = \mathbf{Pr}(\text{Po}(\lambda n t_0) < n) = e^{-\Omega(n \ln n)}$$

by (5). For $k = 0, 1, \ldots$, let $E_k$ be the event $A_{J'_k} \wedge (L_{J'_k}(j^* - 1) \geq \frac{1}{2}n^{1/2}\ln^3 n)$. Let $E = \bigwedge_{k=0}^{n-1} E_k$. We saw earlier that $\mathbf{Pr}(\overline{A_{t_0}}) = o(1)$. By Lemma 5.2, as before, we have $\mathbf{Pr}(L_{J'_k}(j^* - 1) < \frac{1}{2}n^{1/2}\ln^3 n) = e^{-\Omega(\ln^2 n)}$. Thus,

$$\mathbf{Pr}(\overline{E}) \leq \mathbf{Pr}(\overline{A_{t_0}}) + \mathbf{Pr}(J'_n > t_0) + n e^{-\Omega(\ln^2 n)} = o(1).$$

For $k = 0, 1, \ldots$, let $R_k = L_{J'_k}(j^*)$ and for $k = 1, 2, \ldots$, let $Y_k = R_k - R_{k-1}$, so that

$$R_k = R_0 + \sum_{j=1}^{k} Y_j.$$

Let $p_2 = \ln^6 n/(24n)$, and let $r_1 = \lfloor 2\ln^3 n \rfloor$. On $E_{k-1} \wedge (R_{k-1} < r_1)$,

$$\mathbf{Pr}(Y_k = 1|\phi_{J'_{k-1}}) \geq 2p_2$$

and

$$\mathbf{Pr}(Y_k = -1|\phi_{J'_{k-1}}) \leq p_2,$$

for $n$ sufficiently large. [Here we use $\phi_t$ to denote the $\sigma$-field generated by $(X_s : 0 \leq s \leq t)$.] Also, then $np_2 \geq 2r_1$. Hence, by Lemma 7.2, for each integer $r_0$ with $0 \leq r_0 < r_1$,

$$\mathbf{Pr}(E \wedge (R_k < r_1 \; \forall k \in \{1, \ldots, n\})|R_0 = r_0) \leq e^{-p_2 n/28}.$$

Since $\mathbf{Pr}(\overline{E}) = o(1)$, it follows that a.a.s. $R_k \geq r_1$ for some $k \in \{1, \ldots, n\}$. (If $R_0 = r_1$, then we may replace $r_1$ by $r'_1 = r_1 + 1$ above: if $R_0 \geq r_1 + 1$, then $R_1 \geq r_1$ a.s.) Thus, a.a.s. $L_{J'_k}(j^*) \geq \lfloor 2\ln^3 n \rfloor$ for some $k \in \{1, \ldots, n\}$. Finally, since $J'_n \leq t_0$ a.a.s. as we saw above, we find that a.a.s. $L_t(j^*) \geq \lfloor 2\ln^3 n \rfloor$ for some $t \in [0, t_0]$.

In order to complete the proof of (37), it suffices to show that a.a.s. there will be no "excursions" that cross downwards from $\lfloor 2\ln^3 n \rfloor$ to at most $\ln^3 n$. Let $B$ be the event that there is such a crossing. The only possible start times for such a crossing are departure times during $[0, t_0]$. Recall that $N = \lceil 2\lambda n t_0 \rceil$. Now $|X_0| \leq N$ a.a.s. and we saw in (36) that a.a.s. there are at most $N$ arrivals in $[0, t_0]$. Hence, if $C$ denotes the event that there are more than $2N$ departures during $[0, t_0]$, then $\mathbf{Pr}(C) = o(1)$.

We may use Lemma 7.3 (suitably translated and reversed) to upper bound the probability that any given excursion leads to a crossing. Let $a = \lfloor \ln^3 n \rfloor$. Let $p = 2p_2$ and let $q = p_2$. We apply Lemma 7.3 with $a$, $p$, $q$ and $E_k$ as above. We obtain

$$\mathbf{Pr}(B) \leq \mathbf{Pr}(C) + N 2^{-a} + \mathbf{Pr}\left(\bigwedge_{k=0}^{N-1} \overline{E_k}\right) = o(1).$$

Thus, we have established (37), and, hence, proved that $M \geq j^*$ a.a.s.

We now consider upper bounds on $M$. We shall show that $M \leq j^* + 1$ a.a.s., by showing that $L_{t_0}(j^* + 2) = 0$ a.a.s. For $k = 0, 1, \ldots$, let $F_k$ be the event that at the arrival time $J_k$ there are no more than $2n^{1/2} \ln^3 n$ bins with at least $j^*$ balls. Since $l(j^*) < n^{1/2} \ln^3 n$, Lemma 5.1 yields $\mathbf{Pr}(\overline{F_k}) = e^{-\Omega(\ln^2 n)}$. Consider the ball which arrives at time $J_k$: on $F_{k-1}$ it has probability at most $p_3 = 4 \ln^6 n / n$ of falling into a bin with at least $j^*$ balls. Thus, for each positive integer $r$,

$$\mathbf{Pr}(Z_{t_0}(j^* + 1) \geq r) \leq \mathbf{Pr}(B(N, p_3) \geq r) + \mathbf{Pr}\left(\bigvee_{k=0}^{N-1} \overline{F_k}\right) + \mathbf{Pr}(J_{N+1} \leq t_0).$$

Also, the probability that some "initial" ball survives to time $t_0/2$ is at most $\lambda n e^{-t_0/2}$. Hence, there is a constant $c$ such that, with probability $1 - O(n^{-K-3})$, we have $L_t(j^* + 1) \leq c \ln^7 n$ uniformly for all $t \in [t_0/2, t_0]$. Thus, this also holds over $[0, t_0]$.

For $k = 0, 1, \ldots$, let $F_k'$ be the event that at time $J_k$ there are no more than $c \ln^7 n$ bins with at least $j^* + 1$ balls. On $F_{k-1}'$, the ball arriving at time $J_k$ has probability at most $p_4 = c^2 \ln^{14} n \, n^{-2}$ of falling into a bin with at least $j^* + 1$ balls. Then for each positive integer $r$,

$$\mathbf{Pr}(Z_{t_0}(j^* + 2) \geq r) \leq \mathbf{Pr}(B(N, p_4) \geq r) + \mathbf{Pr}\left(\bigvee_{k=0}^{N-1} \overline{F_k'}\right) + \mathbf{Pr}(J_{N+1} \leq t_0).$$

Also, as we noted above, the probability that some "initial" ball survives to time $t_0$ is at most $\lambda n e^{-t_0}$, and so

$$\mathbf{Pr}(M_{t_0} \geq j^* + r + 1) \leq \mathbf{Pr}(Z_{t_0}(j^* + 2) \geq r) + \lambda n e^{-t_0}.$$

It follows on taking $r = 1$ that a.a.s. $M_{t_0} \leq j^* + 1$; and on taking $r = K + 3$ that

$$\mathbf{Pr}(M_{t_0} \geq j^* + K + 4) = o(n^{-K-2}).$$

Now (9) in Lemma 2.1(b), say with $\tau = n^K$ and $s = n^{-2}$, yields the upper bound part of (2).

8.3. *Completing the proof.* In this section $d$ will be any fixed integer at least 2. The lower bound half of (2) will follow from (3), which we now prove. [See also (35) above.] Let $0 < \epsilon < \frac{1}{3}$, and let $\tau = \exp(n^{1/3-\epsilon})$. By Lemma 6.1, there is a constant integer $c > 0$ such that $l(j^* - c) \geq n^{1-\epsilon/2}$. By the concentration result Lemma 5.1 [applied to the function $L(j^* - c)$, with $u = n^{1-\epsilon/2}$],

$$\mathbf{Pr}\,(M < j^* - c) = \mathbf{Pr}\,(L(j^* - c) = 0) = \exp(-\Omega(n^{1/3-\epsilon/3})).$$

Now we may use inequality (8) in Lemma 2.1(b), with $s = 1/\tau$, $a = j^* - c - 3$ and $b = 2$, to show that a.a.s. $M_t \geq j^* - c - 2$ for all $t \in [0, \tau]$. This completes the proof of (3).

It remains to prove (4). Let $z = z(n)$ be a positive integer such that $\ln z = o(\ln n)$. Note that balls choosing bin 1 on each of their $d$ trials arrive in a Poisson process at rate $\lambda n^{-(d-1)}$ (recall that balls choose bins with replacement). Let $C_t$ be the event that, in the interval $[t, t+1)$, there are at least $z$ balls which arrive, choose bin 1 each time, and survive at least to time $t + 1$. Then

$$\begin{aligned}
\mathbf{Pr}\,(C_t) &\geq (1 + o(1))(\lambda n^{-(d-1)} z^{-1})^z e^{-z} \\
&= \exp(-(d-1)z \ln n - z \ln z + O(z)) \\
&= n^{-(d-1+o(1))z}.
\end{aligned}$$

Hence,

$$\begin{aligned}
\mathbf{Pr}\,(M_t^{(n)} < z \,\forall t \in [0, \tau)) &\leq \mathbf{Pr}\,(X_t(1) < z \,\forall t \in [0, \tau)) \\
&\leq \mathbf{Pr}\,(\text{each of } C_0, \ldots, C_{\lfloor \tau \rfloor - 1} \text{ fail}) \\
&\leq (1 - n^{-(d-1+o(1))z})^\tau \\
&\leq \exp(-\tau n^{-(d-1+o(1))z}).
\end{aligned}$$

Hence,

(38) $\mathbf{Pr}\,(M_t^{(n)} < z \,\forall t \in [0, \tau)) \to 0 \quad \text{as } n \to \infty \text{ if } \tau n^{-(d-1+o(1))z} \to \infty.$

This yields (4).

**9. Chaoticity.** As usual, fix a positive integer $d$: let us assume here that $d \geq 2$. One consequence of our concentration results is that asymptotically, as $n \to \infty$, individual bin loads become independent of one another. Thus, our network satisfies the *chaos hypothesis*, Boltzmann's *stosszahlansatz* [6]. In recent years chaoticity phenomena have received considerable attention [6, 9, 10, 18] in the context of various multitype particle systems, such as computer and communication networks, and interacting physical and chemical processes. Consider the equilibrium case.

PROPOSITION 9.1. *Fix an integer $r \geq 2$. Consider the $n$-bin model, with load vector $X$ in equilibrium. For any distinct indices $j_1, \ldots, j_r$, the joint law of $X(j_1), \ldots, X(j_r)$ differs from the product law by at most $O(n^{-1} \ln^4 n)$ in total variation.*

PROOF. As before, let $u(k, X)$ denote the fraction of bins with load at least $k$. By Lemma 5.2,
$$\sup_k \mathbf{Pr}\left(|u(k, X) - u(k)| \geq n^{-1/2} \ln^{3/2} n\right) = O(n^{-1}).$$

Hence, for each positive integer $a \leq r$,

$$\sup_{k_1, \ldots, k_a} \mathbf{E}\left[\prod_{s=1}^{a} |u(k_s, X) - u(k_s)|\right] = O(n^{-a/2} \ln^{3a/2} n) + O(n^{-1}), \quad (39)$$

where the supremum is over all $a$-tuples $k_1, \ldots, k_a$ of nonnegative integers (not necessarily distinct). But

$$\mathbf{E}\left[\prod_{s=1}^{r} u(k_s, X)\right] - \prod_{s=1}^{r} u(k_s)$$
$$= \sum_{A \subseteq \{1, \ldots, r\}, |A| \geq 2} \mathbf{E}\left[\prod_{s \in A} (u(k_s, X) - u(k_s))\right] \prod_{s \in \{1, \ldots, r\} \setminus A} u(k_s).$$

Hence, by (39), uniformly over all $r$-tuples $k_1, \ldots, k_r$,

$$\left|\mathbf{E}\left[\prod_{s=1}^{r} u(k_s, X)\right] - \prod_{s=1}^{r} u(k_s)\right|$$
$$\leq \sum_{A \subseteq \{1, \ldots, r\}, |A| \geq 2} \mathbf{E}\left[\prod_{s \in A} |u(k_s, X) - u(k_s)|\right] \quad (40)$$
$$= O(n^{-1} \ln^3 n).$$

Now
$$u(k, X) = \frac{1}{n} \sum_{j=1}^{n} \mathbf{1}_{X(j) \geq k}.$$

Thus,
$$\mathbf{E}\left[\prod_{s=1}^{r} u(k_s, X)\right] = n^{-r} \mathbf{E}\left[\prod_{s=1}^{r} \sum_{j=1}^{n} \mathbf{1}_{X(j) \geq k_s}\right]$$
$$= \mathbf{E}\left[\prod_{s=1}^{r} \mathbf{1}_{X(s) \geq k_s}\right] + O(n^{-1})$$

uniformly over all $r$-tuples $k_1,\ldots,k_r$, since when we expand the middle expression, there are $O(n^{r-1})$ terms for which the values of $j$ are not all distinct. Hence, from (40),

$$\sup_{k_1,\ldots,k_r} \left| \mathbf{E}\left[ \prod_{s=1}^{r} \mathbf{1}_{X(s) \geq k_s} \right] - \prod_{s=1}^{r} u(k_s) \right| = O(n^{-1} \ln^3 n).$$

But

$$\mathbf{Pr}\left( \bigwedge_{s=1}^{r} (X(s) = k_s) \right) = \mathbf{E}\left[ \prod_{s=1}^{r} (\mathbf{1}_{X(s) \geq k_s} - \mathbf{1}_{X(s) \geq k_s+1}) \right],$$

which is sum of $2^r$ terms $\pm \mathbf{E}[\prod_{s=1}^{r} \mathbf{1}_{X(s) \geq k'_s}]$, where $k'_s = k_s$ or $k_s + 1$; and $\prod_{s=1}^{r} \mathbf{Pr}(X(s) = k_s)$ is a corresponding sum of terms $\pm \prod_{s=1}^{r} u(k'_s)$. Hence, for any set $j_1,\ldots,j_r$ of distinct bin indices,

$$(41) \quad \sup_{k_1,\ldots,k_r} \left| \mathbf{Pr}\left( \bigwedge_{s=1}^{r} (X(j_s) = k_s) \right) - \prod_{s=1}^{r} \mathbf{Pr}(X(j_s) = k_s) \right| = O(n^{-1} \ln^3 n).$$

But there exists a constant $c > 0$ such that

$$\mathbf{Pr}\left( \max_j X(j) > \ln\ln n / \ln d + c \right) = O(n^{-1}).$$

Hence, for any set $j_1,\ldots,j_r$ of distinct bin indices, the joint law of $X(j_1),\ldots,X(j_r)$ differs from the product law by at most $O(n^{-1} \ln^3 n (\ln\ln n)^r)$ in total variation. $\square$

The last result, together with the rapid mixing result Theorem 1.1, shows that, with suitable initial conditions, bin loads will be nearly independent after a short time. Suppose, for example, that we start with all bins empty or, more generally, with $O(n)$ balls in total, and let $t = t(n) \geq 2\ln n$. Let $j_1,\ldots,j_r$ be fixed distinct indices, where $r \geq 2$. Then if $Y$ has the equilibrium distribution,

$$d_{\mathrm{TV}}(\mathcal{L}(X_t(j_1),\ldots,X_t(j_r)), \mathcal{L}(X_t(j_1)) \otimes \cdots \otimes \mathcal{L}(X_t(j_r)))$$
$$\leq d_{\mathrm{TV}}(\mathcal{L}(Y(j_1),\ldots,Y(j_r)), \mathcal{L}(Y(j_1)) \otimes \cdots \otimes \mathcal{L}(Y(j_r)))$$
$$\quad + (r+1) d_{\mathrm{TV}}(\mathcal{L}(X_t), \mathcal{L}(Y))$$
$$= O(n^{-1} \ln^4 n).$$

**10. Concluding remarks.** We have investigated a natural continuous-time balls-and-bins model with $d$ random choices, which exhibits the "power of two choices" phenomenon. We found that the system converges rapidly to its equilibrium distribution; in equilibrium, the maximum load is a.a.s. concentrated on just two values; when $d = 1$, these values are close to $\ln n / \ln \ln n$; and when $d \geq 2$, they are close to $\ln \ln n / \ln d$, and the maximum load varies little over polynomial length intervals. We make three further remarks:

(a) We have not discussed the next level of detail. For example, for given values of $d \geq 2$ and $\lambda > 0$, let $m(d, \lambda; n)$ denote the median value of the maximum load $M^{(n)}$ in equilibrium. We know that the difference $m(d, \lambda; n) - \ln \ln n / \ln d$ stays bounded as $n \to \infty$, but how does it behave in more detail? How does it depend on $\lambda$?

(b) Our approach can be applied, in a natural way, to the "original" load-balancing problem, where $m \sim cn$ balls are thrown into $n$ bins sequentially, each ball chooses $d$ random bins, and is placed into a least loaded of these bins, see [2, 3, 16]. It is known that, with probability tending to 1 as $n \to \infty$, at the end of the allocation process the maximum load of a bin is $\ln \ln n / \ln d + O(1)$, though it has not been possible to determine the behavior of the $O(1)$ term. We make a step forward here, and see that the maximum load is concentrated on at most two values, as in the processes considered earlier in this paper.

We embed the process in continuous time, and for the $n$-bin case, we assume that balls arrive in a Poisson process of rate $n$. A natural coupling, combined with the bounded differences method, yield concentration of measure for Lipschitz functions. As before, let $(X_t)$ denote the loads process, let $u(i, x)$ be the proportion of bins with load at least $i$ in state $x$, and let $u_t(i) = \mathbf{E}[u(i, X_t)]$. Let $t_0 > 0$ be a fixed time. Then uniformly over $t \in [0, t_0]$ and over $i \in \mathbb{N}$,

$$\frac{du_t(i)}{dt} = \mathbf{E}[u(i-1, X_t)^d] - \mathbf{E}[u(i, X_t)^d]$$
$$= u_t(i-1)^d - u_t(i)^d + O(n^{-1} \ln^2 n).$$

Let $(v(t, i) : i = 0, 1, \ldots)$ solve the system of differential equations

(42) $$\frac{dv(t, i)}{dt} = v(t, i-1)^d - v(t, i)^d$$

subject to $v(t, 0) = 1$ for each $t \geq 0$ and $v(0, i) = 0$ for each $i = 1, 2, \ldots$. Then, using Gronwall's lemma (see, e.g., [8]),

$$\sup_{0 \leq t \leq t_0} \sup_{i \in \mathbb{N}} |u_t(i) - v(t, i)| = O(n^{-1} \ln^2 n).$$

Defining $j^* = j^*(n)$ to be the least positive integer $i$ such that $v(t, i) < 2n^{-1/2} \ln n$, with high probability, the maximum load of a bin when about $tn$ balls have been thrown will equal $j^* - 1$ or $j^*$ when $d \geq 3$, and will equal $j^*$ or $j^* + 1$ when $d = 2$. Note that $j^*$ is defined purely in terms of the solution to the limiting differential equation (42).

(c) Our methods can be adapted to handle the "supermarket" model. In this well-studied queueing model, see, for example, [15, 17, 23], there are $n$ single-server queues, with service times which are independent exponentials with mean 1; customers (balls) arrive in a Poisson stream at rate $\lambda n$, where $0 < \lambda < 1$, and go to a shortest of $d$ randomly chosen queues. In [11] we are able to determine (for the first time) the behavior of the maximum queue length, and, indeed, we obtain similar results to those in the present paper. It is possible also to analyze queues with a number $s = s(n)$ of servers, not just 1 or $\infty$.

**Acknowledgment.** We are grateful to a referee for a very detailed reading of the paper.

DEPARTMENT OF MATHEMATICS
LONDON SCHOOL OF ECONOMICS
HOUGHTON STREET
LONDON WC2A 2AE
UNITED KINGDOM
E-MAIL: malwina@planck.lse.ac.uk

DEPARTMENT OF STATISTICS
1 SOUTH PARKS ROAD
OXFORD OX1 3TG
UNITED KINGDOM
E-MAIL: cmcd@stats.ox.ac.uk